\documentclass{amsart}
\usepackage{latexsym,amsxtra,amscd,ifthen}
\usepackage{amsfonts}
\usepackage{verbatim}
\usepackage{amsmath}
\usepackage{amsthm}
\usepackage{amssymb}

\usepackage{enumerate}

\numberwithin{equation}{subsection}

\theoremstyle{plain}
\newtheorem*{theorem}{Theorem}
\newtheorem*{lemma}{Lemma}
\newtheorem*{proposition}{Proposition}
\newtheorem*{corollary}{Corollary}
\newtheorem*{conjecture}{Conjecture}
\newcommand*{\inth}{\textstyle \int}
\theoremstyle{definition}
\newtheorem*{definition}{Definition}
\newtheorem*{example}{Example}

\newtheorem*{remark}{Remark}
\newtheorem*{remarks}{Remarks}
\newtheorem*{question}{Question}



\newcommand{\XXi}{{\xi}}

\DeclareMathOperator{\gldim}{gldim}
\DeclareMathOperator{\THdim}{tHdim}
\DeclareMathOperator{\Hdim}{Hdim}
\DeclareMathOperator{\Hcodim}{Hcdim}
\DeclareMathOperator{\THcodim}{tHcdim}

\DeclareMathOperator{\Ext}{Ext}
\DeclareMathOperator{\Tor}{Tor}

\DeclareMathOperator{\Aut}{Aut}

\DeclareMathOperator{\injdim}{injdim}

\DeclareMathOperator{\GKdim}{GKdim}

\DeclareMathOperator{\D}{{\sf D}}
\DeclareMathOperator{\Hom}{Hom}
\DeclareMathOperator{\RHom}{RHom}

\DeclareMathOperator{\Mod}{{\sf Mod}}

\newcommand{\cal}{\mathcal}

\newcommand{\C}{\mathbb{C}}

\begin{document}

\title[Noetherian Hopf algebras]
{Dualising complexes and \\
twisted Hochschild (co)homology for
\\ Noetherian Hopf algebras}

\author{K.A. Brown and J.J. Zhang}

\address{Brown: Department of Mathematics,
University of Glasgow,
Glasgow G12 8QW, UK}

\email{kab@maths.gla.ac.uk}

\address{zhang: Department of Mathematics, Box 354350,
University of Washington, Seattle, Washington 98195,
USA}

\email{zhang@math.washington.edu}

\begin{abstract}
We show that many noetherian Hopf algebras $A$ have a rigid
dualising complex $R$ with $R \cong \; ^{\nu}A^1 [d]$. Here, $d$
is the injective dimension of the algebra and $\nu$ is a certain
$k$-algebra automorphism of $A$, unique up to an inner
automorphism. In honour of the finite dimensional theory which is
hereby generalised we call $\nu$ the \emph{Nakayama automorphism}
of $A$. We prove that $\nu = S^2\XXi$, where $S$ is the antipode
of $A$ and $\XXi$ is the left winding automorphism of $A$
determined by the left integral of $A$. The Hochschild homology
and cohomology groups with coefficients in a suitably twisted free
bimodule are shown to be non-zero in the top dimension $d$, when
$A$ is an Artin-Schelter regular noetherian Hopf algebra of global
dimension $d$. (Twisted) Poincar{\'e} duality holds in this
setting, as is deduced from a theorem of Van den Bergh.
Calculating $\nu$ for $A$ using also the opposite coalgebra
structure, we determine a formula for $S^4$ generalising a 1976
formula of Radford for $A$ finite dimensional. Applications of the
results to the cases where $A$ is PI, an enveloping algebra, a
quantum group, a quantised function algebra and a group algebra
are outlined.
\end{abstract}





\maketitle


\setcounter{section}{-1}
\section{Introduction}
\label{xxsec0}

\subsection{}
\label{xxsec0.1}
The starting point for the work described in this paper was the
observation by Hadfield and Kr$\ddot{{\text a}}$hmer \cite{HK1}
that, when calculating the Hochschild cohomology of the quantised
function algebra $\mathcal{O}_q (SL(2))$, twisting the coefficient
bimodule $\mathcal{O}_q (SL(2))$ by a suitable algebra
automorphism $\nu$, and so calculating
$H^{\ast}(\mathcal{O}_q (SL(2)),{^{\nu}\mathcal{O}_q (SL(2))^1}),$
avoids the ``dimension drop'' which occurs if one works with
untwisted coefficients. Our initial aim was to understand why
this happens, and to extend their result to other quantised
function algebras and then to other classes of Hopf algebra $A$.

\subsection{Rigid dualising complexes}
\label{xxsec0.2} It transpires that twisting enters the theory at
the level of the \emph{rigid dualising complex} $R$ of $A$,
defined in paragraph (\ref{xxsec4.2}). Let $A$ be a noetherian
Hopf algebra over the base field $k$. (We will assume except where
stated otherwise that $k$ is algebraically closed.) Conjecturally,
$A$ is \emph{Artin-Schelter} (AS-)\emph{Gorenstein} (Definition
1.2): this is true for numerous classes of Hopf algebras, as we
review in $\S$6; see also (\ref{xxsec0.5}). We prove, as
Proposition 4.5:

\begin{theorem}
Let $A$ be a noetherian AS-Gorenstein Hopf $k$-algebra with
bijective antipode. Let $d$ be the injective dimension of $A$.
Then $A$ has a rigid dualising complex $R$,
$$R \cong {^{\nu}}A^1[d].$$
\end{theorem}

Here, ${^{\nu}}A^1$ denotes the $A$-bimodule $A$ with left action
twisted by a certain $k$-algebra automorphism $\nu$ of $A$, which,
following the classical theory of Frobenius algebras which we are
here generalising, we call the \emph{Nakayama automorphism} of
$A$. The automorphism $\nu$ is uniquely determined up to an inner
automorphism.

\subsection{Nakayama automorphisms}
\label{xxsec0.3} To describe the Nakayama automorphism of $A$ we
need to recall from \cite{LWZ} the definition of the \emph{left
integral} of $A$. When $A$ is AS-Gorenstein of dimension $d$,
$\mathrm{Ext}^i_A({_{A}k},{_{A}A})$ is 0 except when $i=d$, and
$\mathrm{Ext}^d_A({_{A}k},{_{A}A})$, the left integral $\inth^l_A$
of $A$, has $k$-dimension 1. Thus
$\mathrm{Ext}^d_A({_{A}k},{_{A}A}) \cong \, {^1k^{\XXi}}$, where
the right $A-$action on the trivial module $k$ is twisted by a
certain left winding automorphism $\XXi$ of $A$. (See $\S\S$1.3
and 2.5.) Let $S$ be the antipode of $A$. We prove, also in
Proposition 4.5:

\begin{theorem} Let $A$ be as in Theorem \ref{xxsec0.2}. Then the
Nakayama automorphism $\nu$ of $A$ is $S^2 \XXi$.
\end{theorem}

\subsection{Twisted Hochschild groups}
\label{xxsec0.4}
Hochschild (co)homology with coefficients in a twisted bimodule
arises naturally when one studies the twisted (co)homology
introduced by Kustermans, Murphy and Tuset \cite{KMT}. Thus, let
$B$ be an algebra of finite global homological dimension, let
$\sigma$ be an algebra automorphism of $B$ and let $HH^\sigma_i(B)$
be the $i$th $\sigma$-twisted Hochschild homology of $B$ (see
\cite{KMT} or \cite{HK1} for details). By \cite[Proposition 2.1]{HK1},
when $\sigma$ is diagonalisable, one has
$$HH^\sigma_i(B)\; \cong \; H_i(B,{^\sigma B}),$$
where the right hand side denotes Hochschild homology with
coefficients in $^\sigma B^1$. We define the {\it
twisted Hochschild dimension} of $B$ to be
$$\THdim B\; =\; \max\{i\;|\; H_{i}(B,{^\sigma B})\neq 0
\quad \text{for some automorphism $\sigma$ of $B$}\}.$$ Similarly
twisted Hochschild codimension, denoted by $\THcodim$, is defined
by using twisted Hochschild cohomology; the definition is in
(\ref{xxsec5.1}). The phenomenon described in (\ref{xxsec0.1}) is
a special case of the following theorem, a central result of this
paper.

\begin{theorem}
[Theorems \ref{xxsec3.4} and \ref{xxsec5.3}] 
Let $A$ be a noetherian AS-Gorenstein Hopf algebra of finite global
dimension $d$, with bijective antipode, and let $\nu$ be its 
Nakayama automorphism. 
Then
\begin{enumerate}
\item
$H_d(A,A^\nu ) \cong Z(A) \neq 0$. As a consequence
$\THdim A=d$.
\item
$H^d(A,{^\nu A}) \cong A/[A,A] \neq 0$. As a consequence
$\THcodim A=d$.
\end{enumerate}
\end{theorem}

In fact, there is a Poincar{\'e} duality connecting the
homology and cohomology of the algebras $A$ of the above
theorem. This is a consequence of a result of Van den Bergh
\cite{VdB1}, which applies in the present context
because of Theorem \ref{xxsec0.2}. The following result is
obtained as Corollary \ref{xxsec5.2}.

\begin{corollary}
Let $A$ be as in Theorem \ref{xxsec0.4}. For every $A$-bimodule
$M$ and for all $i$,
$$H^i(A,M)=H_{d-i}(A,{^{\nu}M}).$$
\end{corollary}

\subsection{Classes of examples}
\label{xxsec0.5}
Our results apply to all currently known classes of noetherian
Hopf algebras (with bijective antipode $S$). In $\S$6 we examine
the details of the application, and in particular the nature of
the Nakayama automorphism, for
\begin{itemize}
\item
affine noetherian Hopf algebras satisfying a polynomial
identity, (\ref{xxsec6.2});
\item
enveloping algebras of finite dimensional Lie algebras, (\ref{xxsec6.3});
\item
quantised enveloping algebras, (\ref{xxsec6.4});
\item
quantised function algebras of semisimple groups,
(\ref{xxsec6.5}),(\ref{xxsec6.6});
\item
noetherian group algebras, (\ref{xxsec6.7}).
\end{itemize}
For all these classes we deduce that the twisted homological and
cohomological dimensions equal the global homological dimension
of the algebra.

As we explain in detail in $\S$6, for some of the above
classes - notably the second and third - the existence and
form of the rigid dualising complex was previously known; and
with it, therefore, the nature of what we are now
calling the Nakayama automorphism. In other cases - notably the
first and fourth - more work will be needed to describe $\nu$.
In the ``classical'' cases of enveloping algebras and group
algebras Poincar{\'e} duality at the level of Lie algebra and
group (co)homology dates back respectively to work of Koszul
\cite{Ko} and Hazewinkel \cite{Ha}, and to Bieri \cite{Bi1}, and
can be retrieved from the present work.

\subsection{The antipode}
\label{xxsec0.6}
If $A$ is a noetherian AS-Gorenstein Hopf algebra we can
apply Theorem \ref{xxsec0.2} also to the Hopf algebra
$(A, \Delta^{\sf op}, S^{-1}, \epsilon)$. The two answers thus
obtained for $\nu$ are necessarily equal to within an inner
automorphism of $A$. Consequently we deduce

\begin{theorem}
[Corollary 4.6]
Let $A$ be as in Theorem \ref{xxsec0.2} with antipode $S$.
Then $S^4 = \gamma \circ \phi \circ \XXi^{-1}$, where
$\gamma$ is an inner automorphism and $\phi$ and $\XXi$ are
the right and left winding automorphisms arising from
the left integral of $A$.
\end{theorem}

This result for finite dimensional Hopf algebras was proved by
Radford \cite{Ra} in 1976 - in that case, $\gamma$ can be
explicitly described and $S$ has finite order. While the latter
corollary is no longer valid when the Hopf algebra is
infinite-dimensional, we can deduce the

\begin{corollary}
[Propositions \ref{xxsec4.6} and \ref{xxsec6.2}]
Let $A$ be as in Theorem \ref{xxsec0.4} with antipode $S$.
Let $io(A)$ be the integral order defined in
\cite[Definition 2.2]{LWZ} and let $o(A)$ be the
Nakayama order of $A$ (see Definition \ref{xxsec4.4}(c)).
\begin{enumerate}
\item
The Nakayama order $o(A)$ is equal to either $io(A)$
or $2\; io(A)$.
\item
Suppose $io(A)$ is finite. Then $S^{4 \; io(A)}$ is an
inner automorphism.
\item
Suppose $A$ is a finite module over its centre. Then
the antipode $S$ and the Nakayama automorphism $\nu$
are of finite order up to inner automorphisms.
\end{enumerate}
\end{corollary}

\subsection{}
\label{xxsec0.7}
Some homological and Hopf algebra background is given in
$\S\S$1 and 2 respectively. Twisted Hochschild homology is
recalled and studied in $\S$3, and Theorem \ref{xxsec0.4}(a)
is proved. Rigid dualising complexes are defined in $\S$4,
and Theorems \ref{xxsec0.2}, \ref{xxsec0.3} and \ref{xxsec0.6}
are proved. In $\S$5 we study twisted Hochschild cohomology,
and prove Theorem 0.4(b) and Corollary 0.4.

\bigskip

\section{General preparations}
\label{xxsec1}

\subsection{Standard notation}
\label{xxsec1.1}

Let $k$ be a commutative base field, arbitrary for the moment.
When we are working with quantum groups it will sometimes be
necessary to assume that $k$ contains ${\mathbb C}(q)$. Unless
stated otherwise all vector spaces are over $k$, and an unadorned
$\otimes$ will mean $\otimes_k$. An algebra or a ring always means
a $k$-algebra with associative multiplication $m_A$ and with unit
1. Every algebra homomorphism is a unitary $k$-algebra
homomorphism. All $A$-modules will be by default left modules. Let
$A^{\sf op}$ denote the opposite algebra of $A$ and let $A^e$
denote the enveloping algebra $A\otimes A^{\sf op}$. The category
of left [resp. right]
$A$-modules is denoted by $A$-$\Mod$ [resp.
$A^{\sf op}$-$\Mod$]. So an $A^e$-module - that is, an object of
$A^e$-$\Mod$ - is the same as an $A$-bimodule central over $k$.

When $A$ is a Hopf algebra we shall use the symbols $\Delta$,
$\epsilon$ and $S$ respectively for its coproduct, counit and
antipode. The coproduct of $a \in A$ will be denoted by
$\Delta(a)=\sum a_1 \otimes a_2$. For details concerning the above
terminology, see for example \cite[Chapter 1]{Mo}. In many places
we assume that
\begin{eqnarray}
\qquad \textit{the antipode $S$ of the Hopf algebra $A$ is
bijective.}\label{antipode}
\end{eqnarray}
Some of the results in this paper are valid more generally, but
the bijectivity of $S$ is satisfied by all the examples
which will concern us. Recently Skryabin proved the following.

\begin{proposition}
\cite{Sk}
Let $A$ be a noetherian Hopf algebra. If one of the following
hold, then $S$ is bijective.
\begin{enumerate}
\item
$A$ is semiprime.
\item
$A$ is affine PI.
\end{enumerate}
\end{proposition}

Thus (\ref{antipode}) is a reasonable condition and may
be satisfied by all noetherian Hopf algebras.

\subsection{Artin-Schelter Gorenstein algebras}
\label{xxsec1.2}

While our main interest is in Hopf algebras of various sorts,
some of our results apply in fact to any algebra possessing
a natural augmentation $\epsilon$ to the base field $k$,
such as, for example, the connected $\mathbb{N}$-graded
$k$-algebras, where $\epsilon$ is the graded map $A\to
A/A_{\geq 1}(=k)$. Accordingly, we give the following
definitions in this more general context.

\begin{definition}
Let $A$ be a noetherian algebra.
\begin{enumerate}
\item
We shall say $A$ has {\it finite injective dimension} if the
injective dimensions of $_AA$ and $A_A$, $\injdim _AA$ and
$\injdim A_A$, are both finite. In this case these integers
are equal by \cite{Za}, and we write $d$ for the common value.
We say $A$ is {\it regular} if it has finite global dimension,
$\gldim _AA < \infty$. Right global dimension always equals
left global dimension \cite[Exercise 4.1.1]{We}; and, when
finite, the global dimension equals the injective dimension.
\item
Suppose now that $A$ has a fixed augmentation $\epsilon: A\to k$.
Let $k$ also denote the trivial $A$-module $A/\ker \epsilon$.
Then $A$ is {\it Artin-Schelter Gorenstein}, which we usually
abbreviate to {\it AS-Gorenstein}, if
\begin{enumerate}
\item[(AS1)]
$\injdim {_AA}=d<\infty$,
\item[(AS2)]
$\dim_k \Ext^d_A({_Ak},{_AA})=1$ and $\Ext^i_A({_Ak},{_AA})=0$
for all $i\neq d$,
\item[(AS3)]
the right $A$-module versions of (AS1,AS2) hold.
\end{enumerate}
\item
If, further, $\gldim A=d$, then $A$ is called {\it Artin-Schelter
regular}, usually shortened to {\it AS-regular}.
\end{enumerate}
\end{definition}

We recall that the following question, first posed in
\cite[1.15]{BG1}, and repeated in \cite{Br1}, remains open:

\begin{question}
Is every noetherian affine Hopf algebra AS-Gorenstein?
\end{question}

We shall note in Section 6 many of the classes of algebra
for which a positive answer to this question is known.

\begin{remark}
The  Artin-Schelter Gorenstein condition (AS2) given in the above
definition is slightly different  from the definition given in
\cite[1.14]{BG1}. When $A$ is a Hopf algebra, the two definitions
are equivalent, as we record in Lemma \ref{xxsec3.2}.
\end{remark}

\subsection{Homological integrals}
\label{xxsec1.3}

Here is the natural extension to augmented algebras of a
definition recently given in \cite{LWZ} for Hopf algebras.
This definition generalises a familiar concept from the
case of a finite dimensional Hopf algebra
\cite[Definition 2.1.1]{Mo}). The motive for the use of the
term ``integral" in the latter setting arises from the
relation to the Haar integral.

\begin{definition}
Let $A$ be a noetherian algebra with a fixed augmentation
$\epsilon: A\to k$. Suppose $A$ is AS-Gorenstein of injective
dimension $d$. Any nonzero element in $\Ext^d_A({_Ak},{_AA})$
is called a {\it left homological integral} of $A$. We denote
$\Ext^d_A({_Ak},{_AA})$ by $\inth^l_A$. Any nonzero element
in $\Ext^d_{A^{\sf op}}({k_A},{A_A})$ is called a {\it right
homological integral} of $A$. We write $\inth^{r}_A=
\Ext^d_{A^{\sf op}}({k_A},{A_A})$. Abusing language slightly, we
shall also call $\inth^l_A$ and $\inth^r_A$ the left and the right
homological integrals of $A$ respectively. When no confusion as to
the algebra in question seems likely, we'll simply write $\inth^l$
and $\inth^r$ respectively.
\end{definition}

\section{Hopf algebra Preparations}
\label{xxsec2}

\subsection{Restriction via $\Delta$}
\label{xxsec2.1}

Let $A$ be a Hopf algebra. Let $\Delta: A\to A\otimes A$ be the
coproduct. Then we may view $A$ as a subalgebra of the algebra
$A\otimes A$ via $\Delta$. Define a functor
$$\mathrm{Res}_A\; : \; (A\otimes A)-\Mod \; \to\;  A-\Mod$$
by restriction of the scalars via $\Delta$, so that
$\mathrm{Res}_A$ is equivalent to the functor $\Hom_{A\otimes
A}({_{A\otimes A}(A\otimes A)}_A,-)$.

\begin{lemma}
 Let $A$ be a Hopf algebra and $B$
be any algebra.
\begin{enumerate}
\item
$\mathrm{Res}_A$ is an exact functor.
\item
Let $N$ be an $A\otimes B^{\sf op}$-module. Then
$\mathrm{Res}_A(A\otimes N)$ with restriction of
$B^{\sf op}$-module from $N$ is isomorphic to
${_AA}\otimes {_k N_B}$ as $ A\otimes B^{\sf op}$-modules
where the left $A$-module on ${_kN_B}$ is trivial.
The isomorphism from $\mathrm{Res}_A(A\otimes N)$ to
${_AA}\otimes {_k N_B}$ is given by
$$a\otimes n\mapsto \sum a_1\otimes S(a_2) n,$$
where $\Delta(a)=\sum a_1\otimes a_2$, with inverse $a'\otimes
n'\mapsto \sum a'_1\otimes a'_2n'$.
\item
$\mathrm{Res}_A$
preserves freeness. Similarly, $A\otimes A$ is a free right
$A$-module via the coproduct $\Delta$.
\item
$\mathrm{Res}_A$
preserves projectivity.
\item
$\mathrm{Res}_A$ preserves
injectivity.
\end{enumerate}
\end{lemma}

\begin{proof}
(a) This is clear.

(b) By the left-module version of the fundamental theorem
of Hopf modules \cite[Theorem 1.9.4]{Mo} $A\otimes N$ is
free over $A$ with basis given by a $k$-basis of $N$,
because $N$ is the space of coinvariants of the left
$A$-comodule $A \otimes N$. It's easy to verify that the
isomorphism given in the proof of \cite[Theorem 1.9.4]{Mo}
has the desired form, and that it preserves the
$B^{\sf op}$-module action.

(c) This follows from (b) and its right-hand analogue, with
$N = A$.

(d) This is true because
$\mathrm{Res}_A(A\otimes A)$ is free by (c).

(e) Let $M$ be an injective $A\otimes A$-module. Then
$$\Hom_A(-,\mathrm{Res}_A(M))\; \cong \; \Hom_A(-,\Hom_{A\otimes A}
({_{A\otimes A}(A\otimes A)_A},M))\; =: \;  (*)$$
By the $\Hom-\otimes$ adjunction
$$(*) \; \cong \; \Hom_{A\otimes A}({_{A\otimes A}(A\otimes A)_A}
\otimes_A -,M)\;  =: \; (**).$$
By (c) $(A\otimes A)_A$ is free, hence flat. Since
$M$ is $A\otimes A$-injective, the functor $(**)$ is exact.
Hence $\mathrm{Res}_A(M)$ is $A$-injective.
\end{proof}

\subsection{The left adjoint action}
\label{xxsec2.2}

The algebra homomorphism $\operatorname{id}_A\otimes S: A\otimes
A\to A\otimes A^{\sf op}=A^e$ induces a functor
$$F_{\operatorname{id}_A\otimes S}\; : \; A^e - \Mod \; \to \;
(A\otimes A) - \Mod .$$ If $S$ is bijective, then 
$F_{\operatorname{id}_A\otimes S}$ is an invertible functor. The
\emph{left adjoint functor} $L$ \cite[Definition 3.4.1(1)]{Mo} is 
defined to be
$$L = \mathrm{Res}_A\circ F_{\operatorname{id}_A\otimes S}\; :\;
A^e-\Mod\; \to \; A-\Mod.$$ Let $M$ be an $A$-bimodule. Then
$L(M)$ is a left $A$-module defined by the action
$$ a\cdot m\; = \; \sum a_1 m S(a_2).$$
 In a similar way to
Lemma \ref{xxsec2.1}, one can prove:

\begin{lemma}
Let $A$ be a Hopf algebra and $L$ be the left adjoint functor
defined above.
\begin{enumerate}
\item
$L$ is an exact functor.
\item
$L$ preserves projectives.
\item
$L(A^e)$ is a free $A$-module.
\item
If $S$ is bijective, then $L$
preserves injectives.
\end{enumerate}
\end{lemma}

\begin{remark}
If $S$ is not bijective, it is unclear to us whether $L$ preserves
injectives. That's why we need bijectivity of $S$ for Lemma
\ref{xxsec2.4}(d).
\end{remark}

\subsection{Twisted (bi)modules}
\label{xxsec2.3}

We extend slightly a standard notation for twisted one-sided
modules, as follows. Let $A$ be an algebra and let $M$ be
an $A$-bimodule. For every pair of algebra automorphisms
$\sigma, \tau$ of $A$, we write $^\sigma M^\tau$ for the
$A$-bimodule defined by
$$a\cdot m\cdot b\; =\; \sigma(a)m\tau(b)$$
for all $a,b\in A$ and all $m\in M$. When one or the other of
$\sigma, \tau$ is the identity map we shall simply omit it,
writing for example $^\sigma M$ for $^\sigma M^1.$ If $\phi$ is
another automorphism of $A$, then the map $x\mapsto \phi(x)$ for
all $x\in A$ defines an isomorphism of $A$-bimodules ${^\sigma
A^\tau}\to {^{\phi\sigma}A^{\phi\tau}}$. In particular,
\begin{eqnarray}{^\sigma A^\tau}\; \cong
\; {A^{\sigma^{-1}\tau}} \; \cong \; {^{\tau^{-1}\sigma}A}
\end{eqnarray}
as $A$-bimodules. If $\tau$ is an inner automorphism of $A$,
given by conjugation $x \mapsto uxu^{-1}$ by the unit
$u$ of $A$, then the map on $A$ given by left multiplication
by $u$ shows that
\begin{eqnarray}
^{\sigma}A^{\beta} \; \cong \; {^{\tau
\sigma}A^{\beta}} \; \cong \; {^{\sigma} A^{\beta\tau}}
\end{eqnarray}
for all automorphisms $\sigma$ and $\beta$ of $A.$

Let $M$ be an $A$-bimodule. Define
$$Z_M (A)\; := \;\{a \in A: am = ma \quad \forall \; m \in M \}.$$
Then $Z_M(A)$ is a subalgebra of $A$. For any  algebra
automorphism $\sigma$, write
$$Z^{\sigma}_A(M)\; =\; \{m\in M\;|\; am=m\sigma(a) \quad
\forall \; a\in A\},$$ so $Z^\sigma (M)$ is a $Z_M(A)$-submodule
of $M$. When $M = A$, $Z^\sigma (M)$ is the space
of $\sigma^{-1}$-{\it normal} elements of $A$. We write $Z(M)$
for $Z^{\operatorname{id}_A}(M)$; of course
$Z(A)=Z_A(A)$ is just the centre of $A$, and $Z^\sigma (A)$ is
a $Z(A)$-module. We shall denote by $N(A)$ the
multiplicative subsemigroup of $\Aut_{k{\text{-}}\mathrm{alg}}(A)$
consisting of those automorphisms $\sigma$ such
that $A$ contains a $\sigma$-normal element which is not a zero divisor.

\begin{lemma}
Let $A$ be an algebra and let $M$ be an $A$-bimodule.
\begin{enumerate}
\item
Under the canonical identification of $M$ with $M^{\sigma}$
as $k$-vector spaces,
$$Z^{\sigma}_A(M)\; =\; Z(M^\sigma).$$
\item Let $\sigma$ and $\tau$ be algebra automorphisms of $A$.
Then $\tau$ carries $Z(A^{\sigma\tau})$ to $Z(A^{\tau\sigma})$.
\item Suppose that $A$ is a Hopf algebra. If $\sigma$ is a Hopf
algebra automorphism of $A$, then $^\sigma(L(M))\cong L(^\sigma
M^\sigma)$. In particular,
$$^\sigma(L(^\tau A))\cong L(^{\sigma \tau}A^\sigma)\cong L(^\tau A).$$
\end{enumerate}
\end{lemma}

\begin{proof} (a) is simply the definition.

(b) By (a) we may regard $Z(A^{\sigma\tau})$ and
$Z(A^{\tau\sigma})$ as subsets of $A$. For every $x\in
Z(A^{\sigma\tau})$, $ax=x\sigma\tau(a)$ for all $a\in A$.
Then, for all $a\in A$,
$$a \tau(x)\; =\; \tau(\tau^{-1}(a)x)\; =\;
\tau(x\sigma\tau\tau^{-1}(a))\;=\;
\tau(x)\tau\sigma(a).$$
This means that $\tau: Z(A^{\sigma\tau})\longrightarrow
Z(A^{\tau\sigma})$ is an isomorphism.

(c) The second assertion is a special case of the first, so we
consider the first one only. Let $\ast$, $\cdot$ and
$\bullet$ be the left action of $A$ on $^\sigma(L(M))$,
$L(M)$ and $L(^\sigma M^\sigma)$ respectively. For every
$a\in A$ and $m\in M$,
$$a\ast m=\sigma(a) \cdot m=\sum \sigma(a)_1
mS(\sigma(a)_2)=:(*).$$
Let $\cdot_\sigma$ be the action of $A$ on
$^\sigma M^\sigma$. Since $\sigma$ commutes with
product, coproduct and antipode, then we continue
$$(*)=\sum \sigma(a_1)
m\sigma(S(a_2))=\sum a_1\cdot_\sigma m\cdot_\sigma
S(a_2)=a\bullet m.$$
\end{proof}

\subsection{Homological properties of the adjoint action}
\label{xxsec2.4}

For us, the usefulness of $L$ rests in its role linking
homological algebra over $A$ with that over $A^e$, thus
permitting the calculation of Hochschild cohomology. The
key lemma is the following result of Ginzburg and Kumar
\cite[Proposition p.197]{GK}:

\begin{lemma}\cite{GK}
Let $A$ be a Hopf algebra and let $M$ be an $A$-bimodule.
\begin{enumerate}
\item
$\Hom_{A^e}(A,M)=\Hom_A(k,L(M))=Z(M)$.
\item
If $S$ is bijective, $\Ext^i_{A^e}(A,M)=\Ext^i_A(k,L(M))$ for
all $i$. 
\end{enumerate}
\end{lemma}

\begin{proof} (a) This was proved in
\cite[Proposition p.197]{GK}. We offer a proof following
the line of the proof of \cite[Lemma 5.7.2(1)]{Mo}. By
definition,
$$\Hom_{A^e}(A,M)\; =\; \{m\in M\;|\; am=ma, \forall \; a\in A\}
\; =\;  Z(M)$$
and
$$\Hom_{A}(k,L(M))\; =\; \{m\in M\;|\; \sum a_1 mS(a_2)=\epsilon
(a)m, \forall \; a\in A\}.$$ If $\sum a_1 mS(a_2)=\epsilon(a) m$
for all $a\in A$, then
$$am\;=\; \sum a_1 m \epsilon(a_2)\;=\; \sum a_1m S(a_2) a_3
\;=\; \sum \epsilon(a_1) ma_2\;=\; ma.$$ Conversely, if $am=ma$
for all $a\in A$, then
$$\sum a_1 mS(a_2)\;=\; \sum a_1S(a_2) m\;=\; \epsilon(a)m.$$
Thus (a) follows.

(b) Let $I$ be an injective $A^e$-resolution of $M$. Since
$L$ is exact and preserves injective modules by Lemma
\ref{xxsec2.2}(d), $L(I)$ is an injective resolution of $L(M)$.
Therefore (b) follows from (a) by replacing $M$ by
its injective resolution. 
\end{proof}

\begin{remark}
In \cite{GK} part (b) was stated for any Hopf algebra. But it
seems to us that the bijectivity of $S$ is needed. This should
also be compared to Proposition \ref{xxsec3.3} in which case the
bijectivity is not needed.
\end{remark}

\subsection{Representations and winding automorphisms}
\label{xxsec2.5}

Let $G(A^\circ)$ be the group of group-like elements of
the Hopf dual $A^\circ$ of the Hopf algebra $A$; that is,
$G(A^\circ)$ is the set of algebra homomorphisms from $A$
to $k$, which is a group under multiplication in
$A^\circ$, namely, convolution of maps, $\pi\ast \pi':=
m_A (\pi\otimes \pi')\Delta$, with identity element the
counit $\epsilon$.  When $A$ is noetherian and is known to be
AS-Gorenstein, so that the left integral $\inth^l$ of $A$ exists
as explained in (\ref{xxsec1.3}), we denote by $\pi_0 \in
G(A^\circ)$ the canonical map
$$\pi_0\; : \; A\to A/\text{r.ann}(\inth^l).$$
Given $\pi \in G(A^\circ)$, let $\Xi^{\ell}[\pi]$ denote the {\it
left winding automorphism} of $A$ defined by
$$\Xi^{\ell}[\pi](a)\; =\; \sum \pi(a_1) a_2$$
for $a\in A$. Similarly, let $\Xi^r[\pi]$ denote the {\it
right winding automorphism} of $A$ defined by
$$\Xi^r[\pi](a)\;=\; \sum a_1\pi(a_2).$$
Since the right winding automorphism will not be used often we
simplify $\Xi^{\ell}[\pi]$ to $\Xi[\pi]$ and call $\Xi[\pi]$ the
winding automorphism associated to $\pi$. It is easy to see that
the inverse of $\Xi[\pi]$ is
$$(\Xi[\pi])^{-1}\;=\; \Xi[\pi S];$$
moreover, $\Xi[\epsilon] =\operatorname{id}_A$, and the map $\Xi:
\pi \longrightarrow \Xi[\pi]$ is an antihomomorphism of groups
$G(A^\circ) \longrightarrow \Aut_{k\text{-}\mathrm{alg}}(A)$. When
$\pi_0$ exists, we shall write
\begin{eqnarray}
\XXi \;:=\; \Xi[\pi_0].
\end{eqnarray}

Given an algebra automorphism (or even an algebra endomorphism)
$\sigma$ of $A$, we define
$\Pi[\sigma]$ to be the map $\epsilon \sigma:A\to k$.
Since $\epsilon  S=\epsilon$,
$$\Pi[S^2]\;=\; \Pi[\operatorname{id}_A]\;=\; \epsilon.$$
Since, in general, $S^2\neq \operatorname{id}_A,$ it follows that
in general $\Xi[\Pi[\sigma]]\neq \sigma$.

Here is a collection of elementary facts about the notations
introduced above.

\begin{lemma}
Let $A$ be a Hopf algebra, let $\pi,\phi \in G(A^\circ)$ and let
$\sigma \in \Aut_{k\text{-}\mathrm{alg}}(A)$.
\begin{enumerate}
\item
$\Pi[\Xi[\pi]]=\pi$. In particular, $\Pi$ is a surjective
map from $\Aut_{k\text{-}\mathrm{alg}}(A)$ to $G(A^\circ)$.
\item
$\Xi[\Pi[\Xi[\pi]]]=\Xi[\pi]$.
\item $\pi  S^2=\pi$ and so $\Xi[\pi S^2]=\Xi[\pi]$. \item
$\Xi[\pi] S^2=S^2  \Xi[\pi]$. In particular, if the left integral
of $A$ exists then $\XXi S^2=S^2 \XXi$.
\end{enumerate}
\end{lemma}

\begin{proof} These are straightforward calculations, so
the proofs are omitted.
\end{proof}

\subsection{The right adjoint action}
\label{xxsec2.6}

Let $M$ be an $A$-bimodule. The {\it right adjoint}
$A$-module $R(M)$ \cite[Definition 3.4.1(2)]{Mo}, which we
shall often denote by $M'$ (see also \cite{FT,HK1}), is
$M$ as a $k$-vector space, with right $A$-action given by
$$m \cdot a\; =\; \sum S(a_2)ma_1$$
for all $a\in A$ and $m\in M$. Below we shall want to
combine twisting with the right adjoint action, producing
right modules of the form $(^\sigma A^\tau)'$ which will
be crucial in our calculation of Hochschild homology.
Parallel results to those stated in (\ref{xxsec2.2}),
(\ref{xxsec2.3}) and (\ref{xxsec2.4}) can be derived for
$R(-).$ We shall omit most of the statements. The following
lemma will be used later. For $\pi \in G(A^\circ)$, we write
$k_{\pi}$ for the corresponding right $A$-module $A/\ker \pi$.

\begin{lemma}
Let $A$ be a Hopf algebra, let $\pi, \phi \in G(A^\circ)$ and let
$\sigma \in \Aut_{k\text{-}\mathrm{alg}}(A)$. Let $M$ be an
$A$-bimodule. From (b) to (e) we assume $S$ is bijective. In (d)
and (e) we assume that the left integral of $A$ exists. In (b) to
(e) we identify the spaces of homomorphisms with vector subspaces
of $M$ or $A$ as appropriate.
\begin{enumerate}
\item $k_{\pi}^\sigma \cong k_{\pi\sigma}$. \item $\Hom_{A^{\sf
op}}(k,(^{S^{-2}}M^\sigma)')\; =\; Z^{\sigma}_A(M)$. \item
$\Hom_{A^{\sf op}}(k_{\phi \ast \pi}, (M^{\Xi[\phi]})') \; = \;
\Hom_{A^{\sf op}}(k_{\pi},(M)')$. As a consequence,
$$\Hom_{A^{\sf op}}(k_{\phi \ast \pi}, (^{S^{-2}}A^{\sigma\Xi[\phi]})')
\; = \; \Hom_{A^{\sf op}}(k_{\pi},(^{S^{-2}}A^{\sigma})').$$
\item
$\Hom_{A^{\sf op}}(\inth^l,(^{S^{-2}}A^\sigma)') \; = \;
Z(A^{\sigma \XXi^{-1}})$. That is,
$$\Hom_{A^{\sf op}}(\inth^l,(A^\sigma)') \;= \;
Z(A^{S^{-2}\sigma \; \XXi^{-1}})
\; \cong \; Z(A^{S^{-2}\XXi^{-1}\sigma}).$$
\item
$\Hom_{A^{\sf op}}(\inth^l, (^{S^{-2}}A^{\XXi} )') \; = \;
Z(A) \; = \; \Hom_{A^{\sf op}} (\inth^l,(^{\XXi^{-1} S^{-2}}A )')$.
\end{enumerate}
\end{lemma}

\begin{proof} (a) follows from the definitions.

(b) By Lemma \ref{xxsec2.3}(a) we may assume
$\sigma=\operatorname{id}_A$. The right $A$-action on
$(^{S^{-2}}M)'$ is then given by
$$ m \cdot a = \sum S^{-1}(a_2)ma_1$$
for all $m \in M$ and $a \in A$, and so we need to prove that
$$\Hom_{A^{\sf op}}(k,(^{S^{-2}}M)')=Z(M).$$
The proof of Lemma \ref{xxsec2.4}(a) can be modified and
we leave the details to the reader.

(c) We may, in a similar way to the proof of (b), identify both
the spaces of maps in question with subspaces of the $k$-vector
space $M$. The rest of the proof is straightforward.

The second assertion follows by taking $M={^{S^{-2}}A^\sigma}$.

(d) As a right module, $\inth^l$ is by definition isomorphic to
$k_{\pi_0}.$ So from (c) with $\pi=\pi_0$ and $\phi=\pi_0^{-1}$,
together with (b) with $A$ for $M$, we see that
$$ \Hom_{A^{\sf op}}(\inth^l,(^{S^{-2}}A^\sigma)') \quad = \quad
Z^{\sigma}_A(A).$$ The first isomorphism in (d) now follows from
Lemma \ref{xxsec2.3}(a). The second isomorphism in (d) follows
from the first after we observe that $(A^{\sigma})' \cong
(^{S^{-2}}A^{S^{-2}\sigma})'$, which is clear from
(\ref{xxsec2.3}.1).

(e) In view of (\ref{xxsec2.3}.2), this is a special case of (d).
\end{proof}

\section{Twisted Hochschild homology}
\label{xxsec3}

\subsection{Definition}
\label{xxsec3.1}

Let $A$ be an algebra and let $\sigma$ be an algebra automorphism.
A $\sigma$-twisted version of the standard cochain complex used to
define cyclic cohomology is defined in \cite[Section 2]{KMT}, and
used to define {\it $\sigma$-twisted cyclic and Hochschild
(co)homology}; the notation is $HC^{\sigma}_{\ast}(A)$ and
$HH^{\sigma}_{\ast}(A)$ in the case of the homology groups, and
analogously for the cohomology groups. As shown in
\cite[Proposition 2.1]{HK1}, when $\sigma$ is diagonalisable, as
will be the case in the most of the examples of particular
interest to us, one has
\begin{eqnarray} HH^\sigma_n(A)\; \cong \;  H_n(A,{^\sigma A}),
\end{eqnarray}
where the right hand side denotes Hochschild
homology with coefficients in $^\sigma A$.

\begin{definition}
The {\it Hochschild dimension} of $A$ is defined to be
$$\Hdim A=\max\{i\;|\; HH_{i}(A)=H_i(A,A)\neq 0\}.$$
The {\it twisted Hochschild dimension of} $A$ is defined
to be
$$\THdim A \;=\; \max\{i\;|\; H_{i}(A,{^\sigma A})\neq 0
\quad \text{for some automorphism $\sigma$ of $A$}\}.$$
\end{definition}

\begin{remark}
Ideally the twisted Hochschild dimension should be defined
in terms of twisted Hochschild homology $HH^\sigma_n(A)$.
Unfortunately our computation in later sections uses
knowledge of $H_{n}(A,{^\sigma A})$; and we don't know in
general when $\sigma$ will be diagonalisable. Further, it's
not immediately obvious whether $H_{n}(A,{^\sigma A}) \neq 0$
is equivalent to $HH^\sigma_n(A)\neq 0$ when $\sigma$ is not
diagonalisable. The above definition is therefore chosen so
as to finesse this uncertainty.

In Section 6 we will comment on the diagonalisability of the
relevant automorphisms for particular classes of examples.
\end{remark}

Of course, in some cases diagonalisability is automatic. The
following lemma is a standard result of linear algebra.

\begin{lemma} If $\sigma$ has finite order and the
characteristic of $k$ is zero, then $\sigma$ is
diagonalisable.
\end{lemma}

\subsection{Ischebeck's spectral sequence}
\label{xxsec3.2}

We need the following special case of Ischebeck's spectral
sequence \cite[1.8]{Is}. Let $A$ be an algebra, let
$M$ be a left $A$-module and $N$ a right $A$-module. Suppose
that $_AM$ has a resolution by finitely generated
projectives and that $N_A$ has finite injective dimension.
Then  Ischebeck's spectral sequence is a convergent
spectral sequence of vector spaces
\begin{eqnarray}
E^{p,q}_2\; :=\; \Ext^p_{A^{\sf op}}(\Ext^{-q}_{A}(M,A),N)
\; \Longrightarrow \; \Tor^A_{-p-q}(N,M).
\end{eqnarray}
If $N$ is an $A$-bimodule, then this is a sequence of left
$A$-modules. In particular, if $A$ has finite injective
dimension, then we obtain the double Ext-spectral sequence
of left $A$-modules
\begin{eqnarray}
E^{p,q}_2\; :=\; \Ext^p_{A^{\sf op}}(\Ext^{-q}_{A}(M,A),A)
\; \Longrightarrow \; {\mathbb H}^{p+q}(M)
\end{eqnarray}
where ${\mathbb H}^n(M)=\begin{cases} 0& n\neq 0\\
M&n=0\end{cases}$ .

As an immediate application we record here the fact that, when $A$
is a noetherian Hopf algebra the (AS2) condition of Definition
(\ref{xxsec1.2})(b) can be replaced by a weaker form, as follows:

\begin{lemma}
Let $A$ be a noetherian Hopf algebra
of injective dimension $d$. If $\Ext^i_A(k,A)$ and
$\Ext^i_{A^{\sf op}}(k,A)$ are finite dimensional over $k$
for all $i$ and if $\Ext^d_A(k,A)\neq 0$, then $A$ is
AS-Gorenstein.
\end{lemma}

\begin{proof} Suppose that the stated hypotheses on the
Ext-groups hold. Using the double Ext-spectral sequence, or
following the proof of \cite[Theorem 1.13]{BG1}, one shows
that $\Ext^i_A(M,A)=0$ and $\Ext^i_{A^{\sf op}}(N,A)=0$ for
all $i\neq d$ and all finite dimensional left $A$-modules
$M$ and finite dimensional right $A$-modules $N$. By the
proof of \cite[Lemma 4.8]{WZ2} $\dim \Ext^d_A(M,A)=\dim M\cdot
\dim \Ext^d_A(k,A)$ for any finite dimension left $A$-module $M$.
Thus
$$\dim \Ext^d_{A}(\Ext^d_{A^{\sf op}}(k,A),A)\; =\;
\dim \Ext^d_{A^{\sf op}}(k,A) \cdot \dim \Ext^d_{A}(k,A).$$
By the double-Ext spectral sequence (3.2.2),
$\Ext^d_{A}(\Ext^d_{A^{\sf op}}(k,A),A)\cong k$. Hence
$$\dim \Ext^d_{A^{\sf op}}(k,A)\; =\; \dim \Ext^d_{A}(k,A)
\; =\; 1.$$
Thus we have proved (AS2).
\end{proof}

\subsection{The link with integrals}
\label{xxsec3.3}

We need two fundamental observations. The first is due
to Feng and Tsygan \cite[Corollary (2.5)]{FT} (see also
\cite[Proposition 2.3]{HK1}):

\begin{proposition}\cite{FT}
Let $A$ be a Hopf algebra and let $M$ be an
$A$-bimodule. Then, for all $i\geq 0$, there are vector space
isomorphisms
$$H_i(A,M)\cong \Tor_i^{A^e}(A,M)\cong \Tor^A_i(M',k).$$
Consequently, $\THdim A\leq \gldim A$.
\end{proposition}

\begin{proof}
By \cite[Lemma 9.1.3]{We},  $H_i(A,M)\cong \Tor_i^{A^e}(A,M)$. To
see the second isomorphism we use the fact that $A\otimes_{A^e}
M\cong M'\otimes_{A}k$ and that the functor $(-)':A^e$-$\Mod\to
A^{\sf op}$-$\Mod$ preserves projectives and is exact. Finally,
$\THdim A\leq \gldim A$ follows from the isomorphisms.
\end{proof}

The second observation is a version of Poincar{\'e} duality
for AS-regular Hopf algebras, which, as we shall note in
Section 6, is sufficient to retrieve some classical results:

\begin{lemma}
Let $A$ be an augmented AS-regular algebra of global dimension $d$.
Then for any right $A$-module $N$,
$$\Tor_{d-i}^A(N,k)\quad \cong \quad \Ext^{i}_{A^{\sf op}}(\inth^l,N)$$
for all $i$.
\end{lemma}

\begin{proof}
The isomorphism follows from Ischebeck's spectral sequence
(3.2.1), taking $M=k$, together with the Artin-Schelter
condition (AS2).
\end{proof}

Noting (3.2.1), and combining the above proposition and lemma,
we immediately deduce the following partial
Poincar{\'e} duality for the Hochschild homology of AS-regular
Hopf algebras.

\begin{corollary}
Let $A$ be an AS-regular Hopf algebra of
global dimension $d$. Let $M$ be an $A$-bimodule. Then,
for every integer $i$,
$$  H_{d-i}(A,M) \; \cong \; \Ext^i_{A^{\sf op}}
(\inth^l,(M)').$$
As a consequence, if $\sigma$ is an algebra
automorphism of $A$, then, for every $i$,
$$  H_{d-i}(A,{^\sigma A}) \; \cong \;
\Ext^i_{A^{\sf op}} (\inth^l,( ^{\sigma}A)').$$
\end{corollary}

\begin{remark} When $S$ is bijective, this corollary can be
deduced from Van den Bergh's general Poincar{\'e} duality
result Proposition \ref{xxsec5.1} -- see Corollary
\ref{xxsec5.2}. However it is unknown at present whether $S$
is bijective in general, so it seems worth recording
this corollary.
\end{remark}

\subsection{Proof of Theorem 0.4(a)}
\label{xxsec3.4}

In this subsection we assume that $S$ is bijective. Recall from
(2.5.1) that $\XXi=\Xi[\pi_0],$ the left winding automorphism
obtained from $\pi_0: A\to A/(\text{r.ann}(\int^l))$.

\begin{theorem}
Let $A$ be a noetherian AS-regular Hopf
algebra of global dimension $d$.
\begin{enumerate}
\item
Let $\sigma$ be an algebra automorphism of $A$. Then
$H_d(A,{^\sigma A})\cong Z({^{S^2\XXi \; \sigma} A})$.
Consequently, $H_d(A,{^\sigma A})\neq 0$ if and only
if $S^2\XXi \; \sigma \in N(A)$.
\item
$H_d(A,{^{\XXi^{-1}S^{-2}}A})\cong Z(A)\neq 0$. Consequently,
$\THdim A=d$.
\item
$\Hdim A=d$ if and only if $S^2\XXi \in N(A)$.
\end{enumerate}
\end{theorem}

\begin{proof} (a) By Corollary \ref{xxsec3.3}
$$H_d(A,{^{\sigma}A}) \; \cong \;
\Hom_{A^{\sf op}}(\inth^l,({^\sigma A})')\; =\;
\Hom_{A^{\sf op}}(\inth^l,(A^{\sigma^{-1}})').$$
By Lemmas \ref{xxsec2.6}(d) and \ref{xxsec2.3}(b),
$$\Hom_{A^{\sf op}}(\inth^l,(A^{\sigma^{-1}})')\;\cong \;
Z(A^{S^{-2}\sigma^{-1}\XXi^{-1}})\;\cong \;
Z(A^{\sigma^{-1}\;\XXi^{-1}S^{-2}})\;\cong \;
Z({^{S^2\XXi\; \sigma}A}).$$
The second assertion follows from the definition $N(A)$

(b,c) Let $M$ be an $A$-bimodule. By Proposition \ref{xxsec3.3},
$H_i(A,M)= \Tor^A_i(M',k)=0$ for all integers $i>d$. Hence $\Hdim
A\leq d$ and $\THdim A\leq d$. Therefore the assertions follow
from (a).
\end{proof}

\subsection{Remarks}
\label{xxsec3.5}

(a) When $A$ is an affine commutative Hopf algebra, and $k$ has
characteristic 0, the global dimension hypothesis always holds
\cite[9.2.11 and 9.3.2]{Mo}, so that $A$ is commutative Gorenstein
and therefore AS-regular \cite{Ba}. Moreover the homological
integral is clearly in this case a trivial module, and $S^2$ is
the identity map by \cite[Corollary 1.5.12]{Mo}. So Theorem 3.4(a)
reduces in this case to the statement that $$HH_d (A) \cong
Z(A)=A,$$ which is a consequence of the
Hochschild-Kostant-Rosenberg theorem \cite[Theorem 9.4.7]{We}.

(b) If $A$ has non-central regular normal elements, then there
are algebra automorphisms $\sigma$ other than
$\XXi^{-1}S^{-2}$ such that $H_d(A,{^{\sigma} A}) \neq 0$.

(c) The determination of the automorphism $\XXi$
in particular classes of examples, and its relation to
other ``canonical'' automorphisms of the algebra,
forms part of the content of later sections.

(d) To understand $H_{d-i}(A,{^\sigma A})$ for all $i$, one will
need to understand the minimal injective resolution of
$({^{S^2\XXi\; \sigma} A})'$. This seems to be quite hard in
general.

\section{AS-Gorenstein algebras and dualising complexes}
\label{xxsec4}

\subsection{Dualising Complexes}
\label{xxsec4.1}

The noncommutative version of the dualising complex was introduced
by Yekutieli in \cite{Ye1}, and has now become one of the standard
homological tools of noncommutative ring theory. We review in the
next four paragraphs the small amount of this theory which we
shall need to draw on. Most of the details, including some of the
necessary facts about derived categories, can be found in
\cite{Ye1, VdB2, YZ}. Let $A$ be an algebra and let $\D^{\rm
b}(A$-$\Mod)$ denote the bounded derived category of left
$A$-modules. A complex $X$ of left $A$-modules is called {\it
homologically finite} if $\bigoplus_i H^i(X)$ is a finitely
generated $A$-module.

\begin{definition}
Let $A$ be a noetherian algebra. A complex $R\in \D^{\rm
b}(A^e$-$\Mod)$ is called a {\it dualising complex} over $A$ if it
satisfies the following conditions:
\begin{enumerate}
\item $R$ has finite injective dimension over $A$ and over $A^{\sf
op}$ respectively. \item $R$ is homologically finite over $A$ and
over $A^{\sf op}$ respectively. \item The canonical morphisms
$A\to \RHom_A(R,R)$ and $A\to \RHom_{A^{\sf op}}(R,R)$ are
isomorphisms in $\D (A^e$-$\Mod)$.
\end{enumerate}
If $A$ is $\mathbb Z$-graded, a graded dualising complex is
defined similarly.
\end{definition}

\subsection{Rigid dualising complexes}
\label{xxsec4.2}

Let $R$ be a complex of $A^e$-modules, viewed as a complex of
$A$-bimodules. Let $R^{\sf op}$ denote the ``opposite complex'' of
$R$ which is defined as follows: as a complex of $k$-modules
$R^{\sf op}= R,$ and the left and right $A^{\sf op}$-module actions on
$R^{\sf op}$ are given by
$$a \cdot r \cdot b \quad := \quad bra$$
for all $a, b\in A^{\sf op}$ and $r\in R^{\sf op}(=R)$. 
If $R\in \D(A^e$-$\Mod)$ then $R^{\sf op}\in \D^{\rm
b}((A^{\sf op})^e$-$\Mod)$. The flip map
$$\tau : (A^{\sf op})^e = A^{\sf op} \otimes A \longrightarrow
A \otimes A^{\sf op} = A^e $$ is an algebra isomorphism, so
$(A^{\sf op})^e$ is isomorphic to $A^e$. Hence there is a natural
isomorphism $\D^{\rm b}(A^e$-$\Mod) \cong
\D^{\rm b}((A^{\sf op})^e$-$\Mod)$.


An unfortunate failing of dualising complexes as defined
in (\ref{xxsec4.1}) is their lack of uniqueness. To remedy
this defect Van den Bergh \cite{VdB2} introduced the idea
of a rigid complex:

\begin{definition}
\cite{VdB2}
Let $A$ be a noetherian algebra.  A dualising
complex $R$ over $A$ is called
{\it rigid} if there is an isomorphism
$$R \cong \RHom_{A^e}(A, R \otimes R^{\sf op})$$
in $\D(A^e$-$\Mod)$. Here the left $A^e$-module structure of
$R\otimes R^{\sf op}$ comes from the left $A$-module structure of
$R$ and the left $A^{\sf op}$-module structure of $R^{\sf op}$.
\end{definition}

Earlier, Yekutieli \cite{Ye1} defined the concept of a
\emph{balanced} dualising complex over a graded algebra $A$; a
balanced dualising complex over a connected graded algebra is
rigid, \cite[Proposition 8.2(2)]{VdB2}. If the algebra $A$ has a
rigid dualising complex $R$, then $R$ is unique up to isomorphism,
\cite[Proposition 8.2]{VdB2}.

\subsection{Van den Bergh condition}
\label{xxsec4.3}

\begin{definition}
Suppose that $A$ has finite injective dimension $d$. Then $A$
satisfies the {\it Van den Bergh condition} if
$$\Ext^i_{A^e}(A,A^e)=\begin{cases} 0& i\neq d\\
U& i=d\end{cases}$$ where $U$ is an invertible $A$-bimodule.
\end{definition}

For us, the Van den Bergh condition will constitute a key
hypothesis in deriving Poincar$\acute{\mathrm{e}}$ duality between
Hochschild cohomology and Hochschild homology; see Proposition
\ref{xxsec5.1}. Our definition is motivated by \cite[Proposition
8.4]{VdB2}, which we reformulate as the following

\begin{proposition}
\cite{VdB2} Let $A$ be a noetherian algebra. Then the Van den
Bergh condition holds if and only if $A$ has a rigid dualising
complex $R = V[s],$ where $V$ is invertible and $s \in
\mathbb{Z}.$ In this case $U = V^{-1}$ and $s = d.$
\end{proposition}

\subsection{Rigid Gorenstein algebras and Nakayama automorphisms}
\label{xxsec4.4}

The definition of the Van den Bergh condition leads us naturally
to the following

\begin{definition}
Let $A$ be an algebra with finite
injective dimension $d$.
\begin{enumerate}
\item
$A$ is {\it rigid Gorenstein} if there is an algebra automorphism
$\nu$ such that
$$\Ext^i_{A^e}(A,A^e)= \begin{cases} 0 & i\neq d\\
{^1A^\nu} &i=d\end{cases}$$ as $A$-bimodules. \item The
automorphism $\nu$ is called the {\it Nakayama automorphism} of
$A$. \item The {\it Nakayama order} of $A$, denoted by $o(A)$, is
the smallest positive integer $n$ such that $\nu^n$ is inner, or
$\infty$ if no such $n$ exists.
\end{enumerate}
\end{definition}

By Proposition \ref{xxsec4.3}, $A$ is rigid Gorenstein if and
only ${^\nu A^1}[d]$ is a rigid dualising complex of
$A$. The terminology in (b) generalises standard usage from
the classical theory of Frobenius algebras, as - for
example - in \cite[Section 2.1]{Ya}. For a Frobenius algebra
$A$ has as its rigid dualising complex the $k$-dual
$A^{\ast}$ of $A$ \cite[Proposition 5.9]{Ye3}, and
$A^{\ast} \cong {^{\nu}A^1}$, by \cite[Theorem 2.4.1]{Ya}.

The following observations are clear from the definition and (2.3.2).

\begin{proposition}
\begin{enumerate}
\item
The Nakayama automorphism is determined up to multiplication
by an inner automorphism of $A.$
\item
The Nakayama automorphism acts trivially on $Z(A).$
\end{enumerate}
\end{proposition}

\subsection{AS-Gorenstein Hopf algebras are rigid Gorenstein}
\label{xxsec4.5}

In this subsection we prove the statement in the heading, and
also that the automorphism $\XXi S^2$ is the Nakayama
automorphism of $A$. (Recall that $\XXi$ is defined in (2.5.1).)
Throughout this subsection we assume that $S$ is bijective.

\begin{lemma}
Let $A$ be a Hopf algebra and let $\pi:A\to k$
be an algebra map. Then, as $A$-bimodules,
$$k_\pi\otimes_A L(A^e)\; = \; k^{\Xi[\pi]}\otimes_A L(A^e)
\; \cong \; {^1A^{\Xi[\pi] S^2}}.$$
\end{lemma}

\begin{proof} The lemma follows from a direct computation
after we understand the module structure on $L(A^e)$.

As a $k$-vector space, $A^e$ is isomorphic to $A\otimes A$; and
there are four different $A$-actions on $A^e$ which commute with
each other. To avoid possible confusions we will try to use
different notations to indicate different places of $A$ and
different actions of $A$. Let $B=A$ as a Hopf algebra and $N=A$ as
$A$-bimodule. We denote the four $A$-module structures as follows.
The left $A$-action on $A$ and the right $B$-action on $A$, both
induced by the multiplication of $A$, are denoted by $\ast_1$ and
$\ast_2$ respectively. The left $B$-action on $N$ induced by the
multiplication of $A$ is denoted by $\ast_3$ and another left
$A$-action on $N$ defined by
$$a\ast_4 n\; =\; nS(a).$$
Hence $A\otimes N$ has a left $A\otimes A$-module structure
induced by $\ast_1$ and $\ast_4$ and a $B$-bimodule structure
induced by $\ast_3$ and $\ast_2$. The left $A$-action on $L(A^e)$
is defined via the coproduct $\Delta$ and the $A\otimes A$-action
on $A\otimes N$. The assertion of the lemma is equivalent to the
following statement: as $B$-bimodules,
\begin{eqnarray}
k_\pi\otimes_A \mathrm{Res}_{\Delta}(A\otimes N)\cong
{^1B^{\Xi[\pi] S^2}}.
\end{eqnarray}

Let $J=\ker \pi$. Then $A/J=k_\pi$ as right $A$-module and
$\pi$ induces a natural $B$-bimodule homomorphism
$$f\; :\;  \mathrm{Res}_{\Delta}(A\otimes N) \; \to \;
\mathrm{Res}_{\Delta}(A\otimes N)/J\;
\mathrm{Res}_{\Delta}(A\otimes N)=k_\pi \otimes_A
\mathrm{Res}_{\Delta}(A\otimes N).$$

To prove the assertion (4.5.1) it suffices to show the following:
\begin{eqnarray}
\Xi[\pi]S^2(b)\ast_3(1\otimes 1)-(1\otimes 1)\ast_2 b\quad = \; \;
\qquad\qquad\qquad\qquad \\
\qquad\qquad \qquad\qquad
\nonumber 1\otimes \Xi[\pi]S^2(b)-b\otimes 1\in J\;
\mathrm{Res}_{\Delta}(A\otimes N).\\
\qquad\qquad \qquad\qquad
b\ast_3(1\otimes 1)=1\otimes b \not\in J\;
\mathrm{Res}_{\Delta}(A\otimes N).
\end{eqnarray}
By Lemma \ref{xxsec2.1}(b), there is an isomorphism
$\Phi: \mathrm{Res}_{\Delta}(A\otimes N)\to A\otimes {_kN}$ where
$\Phi(a\otimes n)=\sum a_1\otimes S(a_2)\ast_4 n$. This
isomorphism induces a commutative diagram
$$\begin{CD}
\mathrm{Res}_{\Delta} (A\otimes N)@>\Phi>> A\otimes {_kN}\\
@V fVV       @VVf'V                             \\
k_\pi\otimes_{A}\mathrm{Res}_{\Delta} (A\otimes N)@>\cong >>
A/J \otimes {_kN}
\end{CD}
$$
where $f': A\otimes N\to A\otimes N/J(A\otimes N)=
A/J\otimes N= k\otimes N$ is the canonical map. Now the
conditions in (4.5.2) and (4.5.3) are equivalent to the
following two conditions
\begin{eqnarray}
f'\Phi(1\otimes \Xi[\pi]S^2(b)-b\otimes 1)=0 \:
\textit{  for all  }\: b\in B.\\
f'\Phi(1\otimes b)\neq 0 \: \textit{  for all  }\: b\in B.
\end{eqnarray}
By definition,
$$f'\Phi(b\otimes 1)=f'(\sum b_1\otimes S(b_2)\ast_4 1)
=f'(\sum b_1\otimes S^2(b_2))=
\sum \pi(b_1)\otimes S^2(b_2).$$
By Lemma \ref{xxsec2.5}(d),
$$\sum \pi(b_1)\otimes S^2(b_2)=\sum \pi(S^2(b_1))\otimes S^2(b_2)
=\sum 1\otimes \pi(S^2(b_1)) S^2(b_2).$$
Since $S^2$ commutes with $\Delta$, we have
$$\sum 1\otimes \pi(S^2(b_1)) S^2(b_2)=1\otimes \Xi[\pi]S^2(b)
=f'\Phi(1\otimes \Xi[\pi]S^2(b)).$$ Thus (4.5.4) follows.
To see (4.5.5) we note, for every $0\neq b\in B$,
$$f'\Phi(1\otimes b)=1\otimes b$$
which is a nonzero element in $k\otimes N$.
\end{proof}

\begin{proposition}
Let $A$ be a noetherian AS-Gorenstein Hopf algebra
(with bijective antipode $S$). Let $d$ be the injective
dimension of $A$.
\begin{enumerate}
\item $A$ is rigid Gorenstein with Nakayama automorphism $\XXi
S^2$. \item The rigid dualising complex of $A$ is $\; {^{\XXi
S^2}A^1}[d]$. \item $\inth^{l}_A \;= \, {^1 k^{\XXi S^2}}\;=\;
{^1k^{\XXi}}.$
\end{enumerate}
\end{proposition}

\begin{proof} (a) By Lemma \ref{xxsec2.4}(b), for all $i \geq 0$,
$$\Ext^i_{A^e}(A,A^e) \cong \Ext^i_A(k,L(A^e)).$$
Since $L(A^e)$ is a free left $A$-module by Lemma 2.2(c), we have
$$\Ext^i_A(k,L(A^e))\cong \Ext^i_A(k,A)\otimes_A
L(A^e).$$
By the AS-Gorenstein condition, for every
$i<d$, $\Ext^i_A(k,A)=0$ and hence $\Ext^i_{A^e}(A,A^e)
=0$. The only nonzero term is
$$\Ext^d_{A^e}(A,A^e)\cong \Ext^d_A(k,A)\otimes_A
L(A^e)\cong k^{\XXi}\otimes_A L(A^e)
\cong {^1 A^{\XXi S^2}}$$
where the last isomorphism follows from the above lemma.
Therefore the assertion follows.

(b) This follows from (a) and Proposition \ref{xxsec4.3}.

(c) This is immediate from the definition, since $S^2$
fixes all one-dimensional representations of $A$.
\end{proof}

\begin{remark}
The case of the above proposition where $A$ has finite
$k$-dimension (or, equivalently, $d=0$), is
\cite[Proposition 3.6]{Sc}. Note also the earlier result
of Oberst and Schneider \cite{OS}: a finite dimensional
Hopf algebra is \emph{symmetric} - that is, its Nakayama
automorphism is the identity - if and only if $S^2 = \XXi
= \operatorname{id}_A.$
\end{remark}

\subsection{The antipode}
\label{xxsec4.6} If $A$ is a Hopf algebra (with bijective
antipode $S$ as usual), then so is $\mathbf{A'} := (A,
\Delta^{\sf op}, S^{-1}, \epsilon )$ \cite[Lemma 1.5.11]{Mo}.
Suppose that $A$ is noetherian and AS-Gorenstein.
Then we can apply Proposition \ref{xxsec4.5} to $\mathbf{A'}$.
In particular we can conclude from (a) of
Proposition 4.5 that $A$ has Nakayama automorphism
$$ \nu' \quad = \quad \phi S^{-2}, $$
where $\phi$ is the \emph{right} winding automorphism
(with respect to $\Delta$) of the left integral of $A$. That
is, in the notation of (\ref{xxsec2.5}),
$$ \phi (a) = \sum a_1 \pi_0 (a_2) $$
for all $a \in A.$

However, as we noted in Proposition \ref{xxsec4.4}(a), the
Nakayama automorphism of $A$ is unique up to an inner
automorphism. In view of Proposition \ref{xxsec4.5}(a) we have
therefore proved the

\begin{corollary}
Let $A$ be a noetherian AS-Gorenstein Hopf algebra with
bijective antipode $S$. Then
\begin{eqnarray}\label{anti}
S^4 \quad = \quad \gamma \circ \phi \circ \XXi^{-1}
\end{eqnarray}
where $\XXi$ and $\phi$ are respectively the left and right
winding automorphisms given by the left integral of $A$, and
$\gamma$ is an inner automorphism.
\end{corollary}

Naturally one immediately asks:

\begin{question}
What is the inner automorphism $\gamma$ in Corollary \ref{xxsec4.6}?
\end{question}

The answer is known when $A$ is finite dimensional, thanks to
Radford's 1976 paper \cite{Ra}. Suppose that $A$ has
finite dimension. Then Radford proved in \cite[Proposition 6]{Ra}
a version of Corollary \ref{xxsec4.6} with the
added information that \emph{$\gamma$ is conjugation by the
group-like element of $A$ which is the character of
the right structure on $\inth_{A^{\ast}}^l.$} In particular,
notice from this that - in general - $\gamma$ is not
trivial. It is tempting to suspect that, for noetherian $A$,
the Hopf dual $A^{\circ}$ of $A$ will play an
important role in answering the above question.

Observe that the three maps composed to give $S^4$ in (4.6.1)
commute with each other, so one deduces at once the
main result of \cite{Ra} - namely, that $S$ has finite order
when $A$ is finite dimensional. We can generalise
this result somewhat, as in the following proposition. Recall
that $io(A)$ denotes the \emph{integral order} of
$A$ which is (by definition) the order of $\XXi$ (see
\cite[Definition 2.2]{LWZ}).

\begin{proposition}
Let $A$ be a noetherian AS-Gorenstein Hopf
algebra with bijective antipode $S$.
\begin{enumerate}
\item The automorphisms $\gamma$, $\phi$ and $\XXi$ in the
corollary commute with each other. \item The Nakayama order $o(A)$
is equal to either $io(A)$ or $2\; io(A)$. \item Suppose $io(A)$
is finite. Then $S^{4 \; io(A)}$ is an inner automorphism. \item
If $A$ has only finitely many one-dimensional modules, then
$io(A)$ is finite. In this case, some power of $S$ and of $\nu$ is
inner.
\end{enumerate}
\end{proposition}

\begin{proof} (a) It is easy to see from the definition that
$\phi$ and $\XXi$ commute. By Lemma \ref{xxsec2.5}(d), $S^2$
commutes with $\XXi$. Similarly, $S^2$ commutes with
$\phi$. Since $S^4=\gamma\circ \phi\circ \XXi^{-1}$, $\phi$ and
$\XXi$ commute with $\gamma$.

(b,c) Let $n=io(A)$ and $m=o(A)$. By the definition of $o(A)$,
$\nu^m$ is inner. Thus, by Proposition 4.5 and Lemma 2.5(d),
$S^{2m}\;\XXi^m $ is inner. Since, for any inner automorphism
$\tau$, $\epsilon \tau=\epsilon$ and since $\epsilon
S^{2m}=\epsilon$, we have
$$\epsilon=\epsilon \nu^m=\epsilon S^{2m} \XXi^m=
\epsilon \XXi^m=\pi_0^m.$$ This means that the order $n$ of
$\pi_0$ in $G(A^\circ)$ divides $m$.

On the other hand, $S^4=\gamma\circ \phi\circ \XXi^{-1}$. Then, by (a),
$$S^{4 n}=\gamma^n\circ \phi^n\circ \XXi^{-n}=\gamma^n.$$
The final equality follows from the fact that $n$ is the order
of $\XXi$ and $\phi$. So $S^{4 n}$ is inner. This
implies that
$$\nu^{2 n}=(S^2 \XXi)^{2n}=S^{4n} \XXi^{2n}=\gamma^n$$
which is inner. Hence $m$ divides $2n$. Therefore $m$ is either
$n$ or $2n$.

(d) Let $A_{ab}$ be the Hopf algebra $A/I,$ where $I$ is the Hopf
ideal generated by elements $xy-yx$ for all $x,y\in A$. If $A$ has
finitely many 1-dimensional modules, then $A_{ab}$ is finite
dimensional. By \cite[Lemma 4.5]{LWZ}, $io(A)$ divides $\dim
A_{ab}$, which is finite. The second statement follows from (b,c).
\end{proof}

\begin{remarks} (1) The conclusion of Proposition 4.6(d) is also
valid when $A$ is
AS-regular and PI (meaning that $A$ satisfies a polynomial
identity), as we'll show in (\ref{xxsec6.2}).

(2) In general both $io(A)$ and therefore $o(A)$ can be infinite
for an AS-Gorenstein noetherian Hopf algebra: consider, for
example the enveloping algebra $A$ of the two-dimensional complex
non-abelian Lie algebra, for which $io(A)$ will be shown to be
infinite in Proposition 6.3(c). Similarly, in general no power of
$S$ for such an $A$ is inner, as illustrated, for example, by $A =
\cal{O}_q(SL(2,\mathbb{C})),$ for a generic parameter $q.$ The
only inner automorphism in this case is the identity, by
\cite[Lemma 9.1.14]{Jo}, and so (6.6.1) confirms that $S^2$ has
infinite order in the outer automorphism group.
\end{remarks}

\section{Twisted Hochschild cohomology}
\label{xxsec5}

In this section we will complete the proof of Theorem
\ref{xxsec0.2}(b).

\subsection{Cohomology and duality}
\label{xxsec5.1}

We begin with some definitions parallel to those introduced for
homology in (\ref{xxsec3.1}). Let $A$ be an algebra, and let
$\sigma$ be an algebra automorphism of $A.$ The \emph{twisted
Hochschild cohomology groups} $HH^{\ast}_{\sigma}(A)$ of  $A$ with
respect to $\sigma$ were defined in \cite{KMT} in a way exactly
analogous to the twisted homology groups discussed in
(\ref{xxsec3.1}). When $\sigma$ is diagonalisable the argument of
\cite[Proposition 2.1]{HK1} applies to show that, for all $i \geq
0$,
$$ HH_{\sigma}^i (A) \; \cong \; H^i (A,\, ^\sigma A).$$
Moreover \cite[Corollary 9.1.5]{We} can be invoked to show
that, for all $i \geq 0,$
$$H^i(A,{^\sigma A})\; \cong \; \Ext_{A^e}^i(A,{^\sigma A}). $$

\begin{definition}
\begin{enumerate}
\item[(a)]
The \emph{Hochschild cohomological dimension} of $A$ is
$$\Hcodim A=\max\{i\;|\; HH^i(A)=H^{i}(A,A)\neq 0 \}.$$
\item[(b)]
The {\it twisted Hochschild cohomological dimension} of $A$ is
$$\THcodim A=\max\{i\;|\; H^{i}(A,{^\sigma A})\neq 0
\; \text{for some automorphism $\sigma$ of $A$}\}.$$
\item[(c)]
We say that $A$ has \emph{finite homological
dimension} if there is an integer $d$ such that
$H^i(A,-)=0$ for all $i>d$.
\end{enumerate}
\end{definition}

In the presence of the Van den Bergh condition introduced in
(\ref{xxsec4.3}) there is a beautiful twisted version of
Poincar$\acute{\mathrm{e}}$ duality for $A$. 
An object $X$ in the derived category $\D (A$-$\Mod)$ is called
{\it compact} if the functor $\Hom_{\D (A-\Mod)}(X,-)$ 
commutes with coproducts. It is well-known that $X$ is 
compact if and only $X$ is quasi-isomorphic to a perfect complex 
over $A$ (namely, $X$ is quasi-isomorphic to a bounded 
complex of finitely generated projective $A$-modules). 

\begin{proposition}\cite{VdB1}
Let $A$ be an algebra of finite homological dimension $d$. 
Suppose that the Van den Bergh condition holds, and
that $A$ is a compact object in $\D (A^e$-$\Mod)$. Then 
for every $A$-bimodule $N$,
$$H^i(A,N)=H_{d-i}(A,U\otimes_A N)$$
for all $i$, where $U$ is the invertible $A$-bimodule
$\Ext^d_{A^e}(A,A^e)$.
\end{proposition}

The above result provides further evidence that, at least 
for rigid Gorenstein algebras, twisted (co)homology is
an important and necessary part of the study of Hochschild 
homology and cohomology. 
Note that the hypothesis of $A$ being 
compact in $\D (A^e$-$\Mod)$ does not appear in \cite{VdB1}; 
it has been pointed out to us by Van den Bergh that this is an
oversight.

\begin{corollary}
Let $A$ be a noetherian rigid Gorenstein algebra with Nakayama 
automorphism $\nu$. Suppose that $A$ has finite
homological dimension $d$ and that $A$ is a compact object in 
$\D (A^e$-$\Mod)$.
\begin{enumerate}
\item
$$H_d(A,{A^{\nu}}) \; \cong\;  Z(A) \neq 0,$$
\noindent and so $\THdim A=d$.
\item
$$H^d(A,{^{\nu}A}) \;\cong \; A/[A,A].$$
\noindent
If $A$ has a non-zero module of finite $k$-dimension
then $\THcodim A=d$.
\end{enumerate}
\end{corollary}

\begin{proof} By definition $U\, = \,{^1A^{\nu}}$. Let
$V = U^{-1}\, =\, {^{\nu}A^1}$.

(a) By Proposition \ref{xxsec5.1},
$$H_d(A,{^1A^{\nu}}) \cong
 H^0(A,A) \cong Z(A)\supset k.$$

(b) By Proposition \ref{xxsec5.1},
$$H^d(A,{^{\nu}A^1}) \; \cong \; H_0 (A,U \otimes_A {^{\nu}A^1})
\;\cong \; H_0(A,A) \; = \; A/[A,A].$$

Suppose that $A$ has a non-zero finite dimensional module. We
claim that $1\not\in [A,A]$, so that $A/[A,A]\neq 0$. To see this
we can pass to a finite dimensional factor ring of $A$ and assume
that $A$ is a simple finite dimensional algebra, hence a central
simple algebra. In that case it is well-known that $1\not\in
[A,A]$.

When $A/[A,A]\neq 0$, then $H^d(A,{^{\nu}A^1})\neq 0$ and
hence $\THcodim A=d$.
\end{proof}

\subsection{Poincar{\'e} duality for Hopf algebras}
\label{xxsec5.2}

In parallel to Proposition \ref{xxsec3.3} for Hochschild
dimension, there is an easy upper bound for the twisted
cohomological dimension of a Hopf algebra. Throughout this
subsection we assume $S$ is bijective.


\begin{lemma}
Let $A$ be a Hopf algebra of finite global dimension $d$.
\begin{enumerate}
\item
$A$ has finite homological dimension, bounded by $d$.
This says that $A$ has a resolution of projective 
left $A^e$-modules bounded at position $-d$. 
Consequently, $\THcodim A\leq d$.
\item
If $A$ is noetherian, then $k$ is a compact object in
$\D(A$-$\Mod)$.
\item
If $A$ is noetherian, then $A$ is a compact object in
$\D(A^e$-$\Mod)$.
\end{enumerate}
\end{lemma}

\begin{proof}
(a) By \cite[Corollary 9.1.5]{We}, for all $i \geq 0,$
$$H^i(A,M)\; \cong \; \Ext_{A^e}^i(A,M)$$
for every $A$-bimodule $M$. By Lemma \ref{xxsec2.4}(b),
$$\Ext_{A^e}^i(A,M)\; \cong \; \Ext^i_A(k,L(M))\; =\; 0$$
for all $i>d$.

(b) This is trivial since $A$ is noetherian and $k$ has 
a bounded resolution of finitely generated projective 
left $A$-modules.

(c) By part (b) $k$ is quasi-isomorphic to a prefect complex 
in $\D(A$-$\Mod)$. Hence, for each $i$, $\Ext^i_{A}(k,-)$ 
commutes with inductive direct limits in $A$-$\Mod$. 
Recall that $A$-$\Mod$ is a Grothendieck category and, 
in particular, inductive direct limits are exact 
\cite[Chapter 2]{Po}.
Since $L(-)$ is a composition of a restriction functor with an
invertible functor, it also commutes with inductive 
direct limits. Hence $\Ext^i_A(k,L(-))$ commutes
with inductive direct limits in $A^e$-$\Mod$. By Lemma 
\ref{xxsec2.4}(b), for any $A^e$-module $M$, 
$$\Ext^i_{A^e}(A,M)\cong \Ext^i_{A}(k,L(M))$$
for all $i$. Therefore $\Ext^i_{A^e}(A,-)$ commutes
with inductive direct limits for all $i$. Together 
with part (a) and the claim in the next paragraph we 
show that $A$ is compact in $\D(A^e$-$\Mod)$.

It remains to show the following: 
{\it if $M$ is a left $A^e$-module of finite projective 
dimension such that 
$\Ext^i_{A^e}(M,-)$ commutes with inductive direct 
limits for all $i$, then $M$ has a projective 
resolution such that each term in the resolution is
a finitely generated projective left $A^e$-module.}
Since $\Hom_{A^e}(M,-)$ commutes with inductive direct 
limits, $M$ is finitely presented by \cite[Theorem 3.5.10]{Po}.
So we have a short exact sequence 
$$0\to M^{-1}\to P^0\to M\to 0$$
where $P^0$ is a finitely generated projective left 
$A^e$-module. Clearly $\Ext^i_{A^e}(P^0,-)$ commutes
with inductive direct limits for all $i$. Using the
long exact sequence we obtain that $\Ext^i_{A^e}(M^{-1},-)$ 
commutes with inductive direct limits for all $i$. 
Then induction on the projective dimension of the 
$A^e$-module $M$ shows that $M$ has a resolution by 
finitely generated projective left $A^e$-modules.
\end{proof}

The following corollary is an immediate consequence
of Propositions \ref{xxsec5.1} and \ref{xxsec4.5}
and the above lemma.

\begin{corollary}
Let $A$ be a noetherian AS-regular Hopf algebra of global 
dimension $d$.
Suppose that $S$ is bijective, and let $\XXi$ be the winding 
automorphism induced by the left integral, (2.5.1). Then for every
$A$-bimodule $M$ and for all $i$, $0 \leq i \leq d,$
$$H^i(A,M)=H_{d-i}(A,{^{\XXi^{-1}S^{-2}}M}).$$
\end{corollary}

\subsection{Proof of Theorem \ref{xxsec0.4}(b)}
\label{xxsec5.3}
 Given an $A$-bimodule $M,$ we let $[A,M]$ denote the
$k$-subspace of $M$ spanned by $\{am-ma : a\in A,\,m\in M\}$.
It is easy to see that
$$H_0(A,M) = A\otimes_{A^e} M=M/[A,M].$$
The following result follows from Corollary \ref{xxsec5.2}, 
and incorporates the remainder of Theorem
\ref{xxsec0.4}(b).

\begin{theorem}
Let $A$ be a noetherian AS-regular Hopf algebra of global 
dimension $d$ and having bijective antipode $S.$ 
Define $\nu= \xi S^2.$
\begin{enumerate}
\item
For every algebra automorphism $\sigma$ of $A,$
$$H^d(A,{^\sigma A})\cong A\otimes_{A^e}{^{\nu^{-1}\sigma}A}
\;=\; {^{\nu^{-1}\sigma}A}/[A,{^{\nu^{-1}\sigma}A}].$$ \item In
particular, $$H^d(A,{^{\nu}A}) \cong A/[A,A]\neq 0,$$ and so
$\THcodim A=d$. \item $\Hcodim A=d$ if and only if
$A^\nu/[A,A^\nu]\neq 0$.
\end{enumerate}
\end{theorem}

\section{Examples}
\label{xxsec6}

\subsection{Some general properties}
\label{xxsec6.1}

To apply the results of Sections \ref{xxsec3} and \ref{xxsec5}
to a particular Hopf algebra $A$, we need to check that

(1) $A$ is AS-regular, and

(2) the antipode $S$ of $A$ is bijective.

In all examples in this section the antipode $S$ is bijective,
thanks either to known descriptions of the antipode, or to
Skryabin's result Proposition \ref{xxsec1.1}, so (2) will present
no difficulties. All the examples considered here are
Auslander-Gorenstein, and other than the group rings of
(\ref{xxsec6.7}), all the examples are Cohen-Macaulay. (See for
example \cite{BG1} for the definitions of these concepts.) So
AS-regularity follows from the next lemma.

\begin{lemma} Let $A$ be a noetherian Hopf algebra.
If $A$ is Auslander-Gorenstein (respectively,
Auslander regular) and Cohen-Macaulay, then it is
AS-Gorenstein (respectively, AS-regular).
\end{lemma}

\begin{proof}
Let $d$ be the injective dimension of $A$. Since $\GKdim k=0$, the
Cohen-Macaulay condition forces $\Ext^i_A(k,A)=0$ for all $i\neq
d,$ and then $\Ext^d_A(k,A)\neq 0$ by \cite[Theorem 2.3(2)]{ASZ1}.
Now the Auslander condition and the Cohen-Macaulay conditions
applied to this latter module show that $\Ext^d_A(k,A)$ has
GK-dimension zero; that is, it is finite dimensional over $k$. The
lemma now follows from Lemma \ref{xxsec3.2}.
\end{proof}

Let $\nu$ be the Nakayama automorphism of $A$. By
(\ref{xxsec3.1}.1), to connect twisted Hochschild homologies
$HH^{\nu}_i(A)$ with $H_i (A,{^{\nu} A})$, one has to check that

(3) $\nu$ is diagonalisable.

In general, this condition will fail - see, for instance,
Proposition \ref{xxsec6.3}(d). We will review the
situation for each class of algebras in turn.

Finally, even when $\nu$ is not the identity, we sometimes
nevertheless find that

(4) $\Hdim A=\Hcodim A=\gldim A$.

We will discuss when this is true.

\subsection{Noetherian PI Hopf algebras}
\label{xxsec6.2}

Let $A$ be a noetherian affine PI Hopf algebra, and consider
conditions (1) and (2) of (6.1). By Skryabin's result Proposition
\ref{xxsec1.1}(b), the antipode of $A$ is bijective, so (2) is
satisfied. As for (1), first note that $A$ is
Auslander-Gorenstein, of dimension $d$, say, by \cite[Theorem
0.1]{WZ1}. By \cite[Theorem A]{Br2} there is an irreducible left
$A$-module $W$ for which $\Ext^d_A(W,A) \neq 0$, so that
\cite[Lemma 1.11]{BG1} confirms that $\Ext^d_A(V,A) \neq 0$ for
every left irreducible $A$-module $V$. Hence \cite[Theorem
3.10(iii)(a)]{SZ} ensures that $A$ is (Krull) Cohen-Macaulay.
Since the Krull and GK-dimensions coincide for affine noetherian
PI algebras by \cite[Proposition 13.10.6 and Theorem 6.4.8]{MR},
we can conclude that $A$ is Auslander-Gorenstein and
Cohen-Macaulay. Thus, by Lemma \ref{xxsec6.1}, $A$ is
AS-Gorenstein.

The situation regarding (3) of (6.1) in the PI case is less
straightforward, since it depends on the characteristic of $k.$
The enveloping algebra of the two-dimensional solvable non-abelian
Lie algebra in positive characteristic demonstrates that
diagonalisability of $\XXi$ and thus of $\nu$ will fail in general
for Hopf algebras which are finite modules over their centres.
Nevertheless, in this setting $\nu$ has finite order up to an
inner automorphism, and so Maschke's theorem comes into play in
characteristic 0. Here is the required result:

\begin{proposition}
Let $A$ be a noetherian affine PI Hopf algebra of finite global
dimension.
\begin{enumerate}
\item
The winding automorphism $\XXi$ has finite order.
\item
Some power of the antipode $S$ of $A$ is inner.
\item
The Nakayama automorphism $\nu$ has finite order up to inner
automorphisms.
\end{enumerate}
\end{proposition}

\begin{proof}
The discussion in the opening paragraph of this
subsection shows that $A$ is homologically homogeneous, so
that $A$ is a finite module over its centre by
\cite[Theorem 5.6(iv)]{SZ}.

(a) By \cite[Lemma 5.3(g)]{LWZ}, $io(A)$ is finite. This
is saying that $\XXi$ has finite order.

(b,c) Follow from part (a) and Proposition \ref{xxsec4.6}(b,c).
\end{proof}

Notice that the proof of the proposition works without the
global dimension hypothesis provided $A$ is a finite
module over its centre.

Given (b) of the proposition, it makes sense to ask:

\begin{question} If $A$ is a noetherian affine PI
Hopf algebra, does the antipode of $A$ have finite order?
\end{question}

For regular PI Hopf algebras the equalities (4) of (6.1) are
valid. To see this, we need a lemma. Recall for its
proof that the Nakayama automorphism acts trivially on the
centre, by Proposition \ref{xxsec4.4}(b).

\begin{lemma}
Let $A$ be a noetherian affine PI Hopf
algebra of finite global dimension, and let
$\nu=\XXi S^2$ be its Nakayama automorphism.
\begin{enumerate}
\item
$\nu, \nu^{-1}\in N(A).$
\item
If $\sigma\in N(A)$ and the restriction of $\sigma$
to the centre is the identity, then
$${^\sigma A}/[A,{^\sigma A}]\;\neq \;0.$$
\end{enumerate}
\end{lemma}

\begin{proof}
(a) By \cite[Corollary 1.8]{BG1}, $A$ is a direct sum of
prime rings. The central idempotents corresponding to the
direct summands of $A$ are fixed by $\nu.$ For this reason
we may assume $A$ is a prime ring (and forget about the
coalgebra structure for the rest of the proof).

Since $A$ is a noetherian affine prime PI algebra, the Goldie
quotient ring $Q(A)$ is isomorphic to the central localization
$AC^{-1}$, where $C=Z(A)\setminus\{0\}$. Clearly $Q(A)$ is a
central simple algebra with centre $Z(A)C^{-1}$ \cite[Theorem
13.6.5]{MR}. The automorphism $\nu$ extends to $Q(A)$ naturally,
still denoted by $\nu$. The restriction of $\nu$ to the centre
$Z(Q)$ is the identity, and so, by the Skolem-Noether theorem,
$\nu$ is an inner automorphism. Thus there is an invertible
element $x\in Q(A)$ such that $\nu(b)=xbx^{-1}$ for all $b\in
Q(A)$. Write $x=ms^{-1}$ where $m\in A$ is not a zero divisor and
$s\in Z(A)\setminus\{0\}.$ Then, for every $a\in A$,
$$\nu(a) m\;=\; x a x^{-1} m\;=\; x a s\;=\; xs a\;=\;
ma.$$
Hence $\nu\in N(A)$. The proof works for any automorphism
whose restriction to the centre is the identity,
and in particular for $\nu^{-1}$.

(b) As in (a) we may assume $A$ is prime. Let $Q=Q(A)$.
It suffices to show that there is an elememt $m\in A$ such
that $m\not\in [Q,{^\sigma Q}]$. By the proof of (a),
$\sigma$ is an inner automorphism, say $\sigma(b)=m b
m^{-1}$ for some regular element $m \in A$. For all $a,b\in Q$,
$$\sigma(b)a-ab=mbm^{-1}a-ab\;=\; m(bm^{-1}a-m^{-1}ab)
\;=\; m[b,m^{-1}a].$$ This means that $[Q,{^\sigma Q}]=m[Q,Q]$.
Since $1\not\in [Q,Q]$, we deduce $m\not\in
[Q,{^\sigma Q}]$.
\end{proof}

\begin{theorem}
Let $A$ be a noetherian affine PI Hopf algebra of
GK-dimension $d$.
\begin{enumerate}
\item
$A$ has rigid dualising complex $^{\nu}A^1 [d]$
where $\nu = \XXi S^2$ is the Nakayama automorphism
of $A$.
\item
Suppose that $A$ has finite global dimension. Then
$\Hdim A=\Hcodim A=\gldim A = d$.
\end{enumerate}
\end{theorem}

\begin{proof}
(a) By the opening paragraph of (\ref{xxsec6.2}) $A$
is AS-Gorenstein of dimension $d$, so this is
immediate from Proposition \ref{xxsec4.5}.

(b) This follows from the above lemma and Theorem 
\ref{xxsec3.4}(c) and \ref{xxsec5.3}(c). Note that the latter
applies because $A$ (and hence $A^e$) are finite 
modules over their affine centres, by \cite[Theorem 5.6(iv)]{SZ}.
\end{proof}

\subsection{Enveloping algebras}
\label{xxsec6.3}

In this subsection we drop the hypothesis that $k$ is
algebraically closed.

\begin{proposition}
Let $\mathfrak{g}$ be a Lie algebra over $k$ of
finite dimension $d$, and let $A$ be its universal
enveloping algebra $U(\mathfrak{g}).$
\begin{enumerate}
\item $S^2=\operatorname{id}_A$. \item $A$ is an AS-regular Hopf
algebra with $\gldim A = d.$ \item \cite{Ye2} The rigid dualising
complex of $A$ is $(A \otimes \bigwedge^d \mathfrak{g})[d].$ Thus
the Nakayama automorphism $\nu = \XXi$ of $A$ is the inverse of
the winding automorphism $\Xi[\chi]$, where $\chi$ is the
representation of $\mathfrak{g}$ on $\bigwedge^d \mathfrak{g}.$
That is,
$$\nu(x) \;= \; x + \mathrm{tr}
(\mathrm{ad}_{\mathfrak{g}}(x))$$
for $x \in \mathfrak{g}.$
\item
If $\mathrm{tr}(\mathrm{ad}_{\mathfrak{g}}(x))\neq 0$,
then $\nu$ is not diagonalisable.
\item
If $\mathfrak{g}$ is semisimple, then $\nu$ is the
identity and hence $\Hdim A=\Hcodim A=d$.
\item
$H_d(A,{A^{\nu}})  \cong  Z(A) \neq 0,$ and
$H^d(A,{^{\nu}A}) \cong  A/[A,A],$ so that
$\THcodim A = \THdim A= d.$ The twisted Poincar{\'e}
duality of Corollary \ref{xxsec5.2} holds for $A$.
\end{enumerate}
\end{proposition}

\begin{proof}
(a) Recall that the antipode is given by $S(x) = -x$ for
$x \in \mathfrak{g}$, so $S^2$ the identity.

(b) That the global dimension is $d$ is
\cite[Theorem XIII.8.2]{CE}. It is well-known that $A$ is
Auslander-regular for all $\mathfrak{g}$.
The assertion follows from Lemma \ref{xxsec6.1}.

(c) The description of rigid dualising complex of $A$
is the main result of \cite{Ye2}.

(d) This follows from (c).

(e) If $\mathfrak{g}$ is semisimple, then $\mathrm{tr}
(\mathrm{ad}_{\mathfrak{g}}(x))=0$ for all $x\in {\mathfrak{g}}$.
By (c), $\nu$ is the identity. The assertion follows from
Theorems \ref{xxsec3.4}(c) and \ref{xxsec5.3}(c).

(f) This follows from Corollaries \ref{xxsec5.1} and \ref{xxsec5.2}.
\end{proof}

\begin{remark}
Twisted Poincar{\'e} duality at the level of ordinary Lie
algebra homology and cohomology follows for
$U(\mathfrak{g})$ from Lemma \ref{xxsec3.3}. This retrieves
the classical duality for Lie algebra (co)homology
which dates back to the work of Koszul \cite{Ko} and
Hazewinkel \cite{Ha}.
\end{remark}

\subsection{Quantised enveloping algebras}
\label{xxsec6.4}

For the definition of the quantised enveloping algebra
$U_q(\mathfrak{g})$ of the semisimple Lie algebra
$\mathfrak{g}$, see for example \cite[I.6.3]{BG2} or
\cite[4.3]{Ja}. The results in this case are already known;
the details are as follows.

\begin{proposition}
Let $\ell$ be an odd positive integer greater than 2. Let
$q$ denote either a primitive $\ell$th
root of one in $\C$, or an indeterminate over $\C$; let $k$ be
$\C$ in the first case, or $\C (q)$ in the second case. Let
$\mathfrak{g}$ be a semisimple complex Lie algebra of dimension
$d$, and let $A$ be the corresponding quantised enveloping algebra
$U_q(\mathfrak g)$. (If $\mathfrak{g}$ has a summand of type $G_2$
assume that 3 does not divide $\ell.$)
\begin{enumerate}
\item $S^2$ is inner. \item $A$ is AS-regular of global dimension
$d$. \item The left integral of $A$ is trivial, i.e.,
$\int^l_A \; \cong\; k$ as $A$-bimodule. Thus $\XXi= \nu =
\operatorname{id}_A$. \item $A$ has rigid dualising complex $A[d]$.
\item
$\Hdim A = \Hcodim A = d$ and Poincar{\'e} duality holds.
\end{enumerate}
\end{proposition}

\begin{proof} (a) See \cite[Proposition 6, p. 164]{KS}
or \cite[(4.9)(1)]{Ja}.

(b,c) That $A$ has finite global dimension was proved in
\cite[Proposition 2.2]{BG1}. Then AS-regularity follows from
\cite[Proposition]{Ch}, and her result shows also that $\gldim A$
is $d$ for all values of the parameter $q,$ and that
$\int_A^l\cong k$. The triviality of $\nu$ now follows from (a)
and Proposition \ref{xxsec4.5}.

(d) This is the main theorem of \cite{Ch}. (Of course it follows
also from (c) and Proposition \ref{xxsec4.5}(b).)

(e) This is immediate from Theorems \ref{xxsec3.4}(c) and
\ref{xxsec5.3}(c), and Corollary \ref{xxsec5.2}, and the fact that
$\nu=\operatorname{id}_A$.
\end{proof}

\subsection{Quantised function algebras}
\label{xxsec6.5}

Our hypotheses on $q$, $\ell$ and $k$ in this paragraph will be the
same as in Proposition \ref{xxsec6.4}. Let $G$ be a connected complex
semisimple algebraic group of dimension $d.$ The definition of the
{\it quantised function algebra} $\mathcal{O}_q(G)$ is given in
the case where $q$ is an indeterminate over $\C$ in
\cite[I.7.5]{BG2}, and in the case where $q$ is a primitive
$\ell$th root of 1 in $\C$ in \cite[III.7.1]{BG2}.

\begin{proposition}
Let $G$, $d$, $q$, $\ell$ and $k$ be as stated
above. Let $A = \mathcal{O}_q(G).$
\begin{enumerate}
\item
$A$ is AS-regular of dimension $d.$
\item
$\nu = \XXi S^2$, $A$ has rigid dualising complex
$^{\nu}A^1 [d]$, $\THdim A=\THcodim A=d$, and twisted
Poincar{\'e} duality holds.
\end{enumerate}
\end{proposition}

\begin{proof}
(a) Suppose first that $q$ is a root of unity.
By \cite[Theorem 2.8]{BG1} $A$ is noetherian PI, Auslander
regular and Cohen Macaulay of global dimension $d$, and
everything is proved in (\ref{xxsec6.2}). Suppose now
that $q$ is an indeterminate. It is proved in
\cite[Proposition 2.7]{BG1} that $\gldim A$ is finite,
and indeed the global dimension is $d$ by
\cite[Theorem 0.1]{GZ}. By \cite[Theorem 0.1]{GZ},
$A$ is Auslander-regular and Cohen-Macaulay. By
Lemma \ref{xxsec6.1}, $A$ is AS-regular.

(b) This follows from Proposition \ref{xxsec4.5} and
Corollaries \ref{xxsec5.1} and \ref{xxsec5.2}.
\end{proof}

\subsection{The special linear group}
\label{xxsec6.6}

In Proposition \ref{xxsec6.5} we did not answer the following

\begin{question}
What is the precise form of the Nakayama automorphism
$\nu=\XXi S^2$ of $\mathcal{O}_q(G)$?
\end{question}

Here we answer the question for the case $G = SL(n,\C)$, giving
$\nu$ in a form which, we conjecture, remains valid for \emph{all}
semisimple groups $G.$ First we need the following lemma, a
straightforward exercise in the use of change of rings formulae
proved in \cite[Lemma 2.6]{LWZ}. Recall the definition of
$\tau-$normal elements in (\ref{xxsec2.2}).

\begin{lemma}
Let $A$ be an AS-Gorenstein algebra of injective
dimension $d$, and let $x$ be a $\tau$-normal non-zero-divisor in
the augmentation ideal of $A$. Let $\bar{A}:=A/\langle x \rangle.$
\begin{enumerate}
\item If $M$ is an $x$-torsionfree left $A$-module and $N$ is an
$x$-torsionfree right $A$-module, then
$$\Ext^1_A(A/(x),M)={^{\tau}(M/xM)},
\quad \Ext^1_{A^{\sf op}}(A/(x),N)=(N/Nx)^{\tau^{-1}}.$$
\item
$\bar{A}$ is an AS-Gorenstein algebra of dimension $d-1$,
with augmentation induced by that of $A$.
\item
$\quad \quad \Ext^d_A(k^{\tau^{-1}},A)\quad \cong \quad
(\Ext^d_A(k,A))^{\tau}\quad \cong \quad
\Ext^{d-1}_{\bar{A}}(k,\bar{A}).$
\item
As bimodules, $ \qquad \inth^l_A\quad \cong \quad
(\inth^l_{\bar{A}})^{\tau^{-1}}$ .
\end{enumerate}
\end{lemma}

Recall, from \cite[I.2.2-I.2.4]{BG2} for example,
that $\mathcal{O}_q (SL_n)$ is generated by elements $X_{ij}$
for $i,j = 1, \ldots , n,$ with relations as given there.
(The key ones for the proof of the proposition are
stated at the start of the proof below.)

\begin{proposition}
Let $n$ be an integer greater than 1, and let
$A={\mathcal O}_q(SL_n)$ where $q$ is a nonzero scalar.
\begin{enumerate}
\item
$\inth^l_A \; \cong \; A/\langle X_{ii}-q^{2(n+1-2i)},
X_{tj}: 1 \leq i,j,t \leq n, \; t\neq j \rangle ,$ as right
$A$-modules.
\item
The automorphism $\XXi$ is determined by
$$\XXi(X_{ij})=q^{2(n+1-2i)}X_{ij}$$
for all $i,j$.
\item
The Nakayama automorphism of $A$ is determined by
$$\nu(X_{ij})=q^{2(n+1-i-j)}X_{ij}$$
for all $i,j$. If $q$ is indeterminate,
then $\nu$ is unique.
\item
$\nu$ is diagonalisable.
\end{enumerate}
\end{proposition}

\begin{proof}
(a) Recall that the relations involving $X_{1n}$ and $X_{n1}$ are
$$X_{1j}X_{1n}=qX_{1n}X_{1j}: j\neq n \quad
\text{and}\quad X_{in}X_{1n}=q^{-1}X_{1n}X_{in}: i\neq 1,$$
$$X_{j1}X_{n1}=qX_{n1}X_{j1}: j\neq n \quad
\text{and}\quad X_{ni}X_{n1}=q^{-1}X_{n1}X_{ni}: i\neq 1;$$ and
$$X_{1n}X_{ij}=X_{ij}X_{1n}\quad \text{and}\quad
X_{n1}X_{ij}=X_{ij}X_{n1}: i\neq 1, j\neq n.$$
In particular, therefore, $X_{1n}$ and $X_{n1}$ are normal
elements. Let $\bar{A}=A/(X_{1n})$ and $\hat{A}=\bar{A}/(X_{n1})$.
Consider first $X_{1n}.$ The relations show that this element
is $\tau$-normal, where
$\tau: X_{11}\to q^{-1}X_{11}, X_{ii}\to X_{ii}$ for $i\neq 1,n$,
$X_{nn}\to qX_{nn}$, and, when $i\neq j$,
$X_{ij}\to b_{ij}X_{ij}$ for some non-zero scalars $b_{ij}$
whose exact value need not concern us. By Lemma
\ref{xxsec6.6}(d),
$$ \inth^l_A\quad \cong \quad (\inth^l_{\bar{A}})^{\tau^{-1}}.$$
The argument for $X_{n1}$ is similar: again, one finds that
$$ \inth^l_{\bar{A}}\quad \cong \quad (\inth^l_{\hat{A}})^{\tau^{-1}}.$$
Notice now that in $\hat{A}$, the cosets represented by
$\{X_{2n},X_{n2}\}$ form a pair of normal elements similar
to $\{X_{1n},X_{n1}\}$. Proceeding as above, a similar relation
between integrals holds. The next pair of normal
elements to deal with is $\{X_{1n-1},X_{n-11}\}$. Continuing in this
way to deal with all off-diagonal generators
in pairs, eventually all such $X_{ij}$s are factored out.
Let $\theta$ denote the composition of two copies of the
winding automorphisms associated to each $X_{ij},$ for $i < j,$ and
set $B=k[X_{11}^{\pm
1},\ldots,X_{(n-1)(n-1)}^{\pm 1}].$ We reach the conclusion
$$ \inth^l_A \quad = \quad \Ext^{n^2-1}_A(k,A)
\quad \cong \quad (\Ext^{n-1}_B( k, B))^{\theta^{-1}}.$$
That is, as right $A-$modules,
$$\inth^l_A \quad \cong \quad A/\langle X_{ii}-q^{2(n+1-2i)},
X_{tj}: 1 \leq i,j,t \leq n, \; t\neq j \rangle ,$$
as claimed.

(b) Since left winding automorphisms of $\mathcal{O}_q (SL_n)$
take a constant value on each column of the matrix
$(X_{ij})$, this follows from (a) and the definition of $\XXi$.

(c) This is immediate from (b) and Proposition \ref{xxsec6.5}(b),
given that
\begin{eqnarray} S^2(X_{ij})=
q^{2(i-j)}X_{ij},\end{eqnarray} \cite[Theorem 4]{FRT}. When $q$
is indeterminate, $\mathcal{O}_q (SL_n)$ does not
have any non-trivial units, \cite[Lemma 9.1.14]{Jo}. Hence
$\nu$ is unique.

(d) This is clear by using a monomial basis of
$\mathcal{O}_q (SL_n)$.
\end{proof}

\begin{remarks}
(a) The argument we've just used for Proposition \ref{xxsec6.6}
works also for
${\mathcal O}_q(GL_n)$, and one gets the identical conclusion.

(b) The same method applies to the multiparameter quantum group
${\mathcal O}_{\lambda,p_{ij}}(GL_n)$, defined,
for example, in \cite[I.2.2-I.2.4]{BG2}. We leave the details
to the interested reader.

(c) Hadfield and Kr\"{a}hmer in \cite{HK2} obtain Theorem
\ref{xxsec0.4}(a) for $\mathcal{O}_q(SL_n)$ with the
same automorphism $\nu$ as in Proposition \ref{xxsec6.5}(b).

(d) We give now a more symmetric, coordinate free formula
for the Nakayama automorphism $\nu$ which, we
conjecture, remains valid for arbitrary $G.$ Recall,
from \cite[I.9.21(a)]{BG2} for example, that given a Hopf
algebra $H$ there are left and right module algebra actions
of the Hopf dual $H^{\circ}$ on $H$, defined
respectively by
$$ u.h = \sum h_1 \langle u,h_2 \rangle \quad \mathrm{ and }
\quad h.u = \sum \langle u, h_1  \rangle h_2, $$
for $h \in H$ and $u \in H^{\circ}.$ In particular, the right
and left
winding automorphisms of $H$ occur in this way, when we take
$u$ to be a group-like element of $H^{\circ}.$ When
$H = \mathcal{O}_q(G)$ and $\mathfrak{g}$ is the Lie algebra
of $G$, $U_q(\mathfrak{g}) \subseteq H^{\circ},$ and
the group-like elements contained in $U_q(\mathfrak{g})$ are
precisely the group $\langle K_{\alpha} : \alpha \in
P \rangle,$ where $P$ is the root lattice of $\mathfrak{g}.$

Now let $\rho \in P$ be half the sum of the positive roots,
and consider the right action of $K_{4 \rho}$ on
$\mathcal{O}_q(SL(n,\C)).$
We calculate that, for $i,j = 1, \ldots , n,$
$$ X_{ij}.K_{4 \rho} = \sum_t X_{it}(K_{4 \rho})X_{tj} =
q^{(\omega_i - \omega_{i-1},4 \rho)}X_{ij} = q^{2(n-2i+1)}X_{ij}.$$
Here, $\omega_1, \ldots , \omega_{n-1}$ are
the fundamental weights of $\mathfrak{sl}(n,\C),$ and
$\omega_0 = \omega_n = 0.$ Comparing this with Proposition
\ref{xxsec6.6}(b) and bearing in mind Theorem \ref{xxsec6.5}(b)
and the notation of (2.5), we see that
\begin{eqnarray}
\pi_0 = K_{4 \rho} \quad \textit{  for the case } \quad G = SL(n,\C).
\end{eqnarray}

\begin{conjecture} With the above notation, (6.6.2) is valid
for all semisimple groups; that is, the Nakayama automorphism
of $\mathcal{O}_q(G)$ is $S^2 \Xi[K_{4 \rho}]$ for
all semisimple groups $G$ and all values of $q$ permitted
in Proposition \ref{xxsec6.4}.
\end{conjecture}

\end{remarks}

\subsection{Noetherian group algebras}
\label{xxsec6.7}

\begin{theorem}
Let $k$ be a field of characteristic $p \geq 0,$ and
let $G$ be a polycyclic-by-finite group of Hirsch length $d.$
\begin{enumerate}
\item $kG$ is noetherian Auslander Gorenstein and AS-Gorenstein of
injective dimension $d.$ \item $kG$ is AS-regular if and only if
it is Auslander regular if and only if $G$ contains no elements of
order $p.$
\end{enumerate}
\end{theorem}

\begin{proof} (a) First, $kG$ is noetherian by
\cite[Corollary 10.2.8]{Pa}. Since we can find a
poly-(infinite cyclic) normal subgroup $H$ of finite index
in $G$ \cite[Lemma 10.2.5]{Pa}, it is sufficient to
prove that (i) $kH$ is AS-Gorenstein (respectively
Auslander Gorenstein) when $H$ is poly-(infinite cyclic),
and (ii) $kH$ AS-Gorenstein (respectively
Auslander Gorenstein) implies that
$kG$ is AS-Gorenstein (respectively
Auslander Gorenstein) when $\mid G : H \mid$ is finite.

(i) We prove this by induction on the Hirsch number $d$ of $H.$
When $d = 0$ there is nothing to prove, so assume the result is
true for all poly-(infinite cyclic) groups of Hirsch number less
than $d$. Choose a normal subgroup $T$ of $H$ with $H/T =\langle
xT \rangle$ infinite cyclic.

AS-Gorenstein property: By the induction hypothesis,
$\mathrm{Ext}^i_{kT}(k,kT) = \delta_{i,d-1}k.$ Applying $kH
\otimes_{kT} - $ to a projective $kT-$resolution of $k$, it is
easy to prove that, for all $i \geq 0,$
\begin{eqnarray}
 \mathrm{Ext}^i_{kH}(kH \otimes_{kT} k, kH) \quad
\cong \quad \mathrm{Ext}^i_{kT}(k,kT)\otimes_{kT} kH
 \end{eqnarray}
as right $kH-$modules. In particular, therefore, as
right $k\langle x \rangle-$modules,
$$\mathrm{Ext}^i_{kH}(kH \otimes_{kT} k, kH) \cong
\delta_{i,d-1}k\langle x \rangle.$$
Let's write
$^H_T\negmedspace\uparrow \negmedspace k$ for $kH \otimes_{kT} k.$
The short exact sequence of left $kH-$modules
$$
0 \longrightarrow ^H_T\negmedspace \uparrow \negmedspace k
\longrightarrow ^H_T\negmedspace \uparrow \negmedspace
k \longrightarrow k \longrightarrow 0,
$$
where the left-hand embedding is given by right multiplication
by $x-1,$ induces exact sequences
$$
\mathrm{Ext}^i_{kH}(^H_T\negmedspace \uparrow \negmedspace k, kH)
\longrightarrow
\mathrm{Ext}^i_{kH}(^H_T\negmedspace \uparrow \negmedspace k, kH)
\longrightarrow
\mathrm{Ext}^{i+1}_{kH}(k,kH)\longrightarrow
\mathrm{Ext}^{i+1}_{kH}(^H_T\negmedspace \uparrow \negmedspace k,
kH),$$
where the left hand and right hand maps are left multiplication
by $x-1,$ \cite[Lemma 2.1]{Br2}. Taking $i< d-2$ here, and applying
the induction hypothesis together with (6.7.1), shows that
$\mathrm{Ext}_{kH}^{t}(k,kH) = 0$ for $t \leq d-2.$ Taking
$i = d-2,$ we find in view of (6.7.1) that
$\mathrm{Ext}_{kH}^{d-1}(k,kH)$ is the kernel of the map
$$\mathrm{Ext}^{d-1}_{kH}(^H_T\negmedspace \uparrow \negmedspace
k,kH) \to \mathrm{Ext}^{d-1}_{kH}(^H_T\negmedspace \uparrow
\negmedspace k,kH),$$
a homomorphism $\rho$ between
free $k\langle x \rangle-$modules of rank one given by
multiplication by $x-1.$ Therefore,
$\mathrm{Ext}_{kH}^{d-1}(k,kH)$ is $0.$ Moreover,
since $\Ext^{d+1}_{kH}(^H_T\negmedspace \uparrow
\negmedspace k, kH)$ is $0$, by (6.7.1), the cokernel of
$\rho$ is $\mathrm{Ext}^d_{kH}(k,kH)$. Thus $\mathrm{Ext}^d_{kH}(k,kH)$
has $k$-dimension one, as required.

Auslander Gorenstein property: Note that $kH=kT[t^{\pm 1};\sigma]$
for some automorphism $\sigma$ of $kT$. By the induction
hypothesis, $kT$ is Auslander Gorenstein. By \cite[Theorem 4.2]{Ek},
$kT[t;\sigma]$ is Auslander Gorenstein. By \cite[Proposition 2.1]{ASZ2}
the localization $kT[t^{\pm 1};\sigma]$ is Auslander Gorenstein.

(ii) This follows because crossed products of finite groups
are Frobenius extensions - see \cite[Lemma 5.4]{AB}.

(b) This follows from (a) and \cite[Theorem 10.3.13]{Pa}.
\end{proof}

\begin{remarks}
\begin{enumerate}
\item
It is an open question whether all noetherian group
algebras $kG$ come from polycyclic-by-finite
groups $G$.
\item
In general a noetherian group algebra $kG$ could have
infinite GK-dimension \cite[Corollary 11.15]{KL}.
Even when $\GKdim kG$ is finite, it could be larger
than the global dimension of $kG$. Hence in these cases,
$kG$ is not Cohen-Macaulay.
\end{enumerate}
\end{remarks}

Poincar{\'e} duality of group homology and cohomology was
investigated by Bieri in a series of papers beginning
with \cite{Bi1}. Given a commutative ring $R$, a non-negative
integer $n$ and a group $G$, he defined
\cite[5.1.1]{Bi2} $G$ to be a \emph{duality group of
dimension $n$ over $R$} if there exists an $RG-$module $C$
such that, for each $RG-$module $A$ and each $k \in \mathbb{Z}$,
there are isomorphisms
$$ H^k(G,A) \cong H_{n-k}(G, C \otimes_R A) $$
induced by the cap product. Moreover, \cite[Definition 3.2]{Bi2},
$G$ is a \emph{Poincar{\'e} duality group of
dimension $n$ over $R$} if $C$ above is $R$ with
\emph{almost trivial} $G-$structure, where this means that the
image of $G$ in $\mathrm{Aut}_R(R)$ is $\{\pm \mathrm{Id} \}.$

It's clear that, for polycyclic-by-finite groups, Bieri's
results describe the shadow  at the level of group
(co)homology of the existence and nature of the dualising
complex for $kG$. On the other hand it is
straightforward to apply the results of \cite{Bi1}, \cite{Bi2}
to determine the Nakayama automorphism for $kG$. In
order to do this, we define the \emph{adjoint trace} $\texttt{T}$
of the polycyclic-by-finite group $G$ as
follows. Fix a series
\begin{eqnarray}
 1 = G_0 \subseteq G_1 \subseteq \cdots \subseteq G_n = G
\end{eqnarray}
 of normal subgroups $G_i$ of $G$ whose factors $G_i/G_{i-1}$
are either torsion-free Abelian, or finite. Provide
 each torsion-free factor $A_i = G_i/G_{i-1}$ with a
$\mathbb{Z}-$basis (of cardinality $t_i$, say), so that
conjugation by $g \in G$ assigns
 a matrix $\tilde{g}_i \in \mathrm{GL}_{t_i}(\mathbb{Z})$
to the pair $(g, A_i).$ Now define a map
 $$ \texttt{T}: G \longrightarrow \{\pm 1 \} : g \mapsto
\prod_i \mathrm{det} (\tilde{g}_i) , $$
 where the product is taken over those $i$, $1 \leq i \leq n,$
such that $G_i/G_{i-1}$ is infinite. One checks
 easily that this definition doesn't depend on the choice of
series (6.7.2), and that $\texttt{T}$ is a group
 homomorphism \cite [\S 3.2]{Bi1}. Finally, set $\texttt{T}_k$
to be the composition of $\texttt{T}$ with the
 canonical homomorphism from $\{\pm 1 \}$ to $k$, and use the
same symbol for the resulting representation of $kG$
 on $k.$

\begin{proposition}
Let $k$ be a field of characteristic $p \geq 0,$ and let $G$
be a polycyclic-by-finite group with Hirsch length
$d.$ If $p > 0$ assume that $G$ contains no elements of
order $p.$
\begin{enumerate}
\item
\cite{Bi2}
$G$ is a Poincar{\'e} duality group over $k$. The dualising
module $C$ is $H^d(G,kG);$ that is, $C$ is the left
integral of $kG.$
\item
The Nakayama automorphism $\nu$ of $kG$ is the winding
automorphism induced by the adjoint trace
$\texttt{T}_k$ of $G$ in $k$. That is,
$$ \nu (g) \quad = \quad \texttt{T}_k (g) g $$
for all $g \in G.$ In particular, $\nu$ is diagonalisable.
\item
$H_d(kG,kG^{\nu}) \cong Z(kG) \neq 0,$ and $H^d(kG,\,^{\nu}kG)
\cong kG/[kG,kG] \neq 0.$
\item
The integral order of $kG$ is either 1 or 2.
\end{enumerate}
\end{proposition}

\begin{proof}
(a) The first statement is proved in
\cite[Theorems 5.6.2 and 5.6.4]{Bi2}, and the second
is \cite[Proposition 5.2.1]{Bi2}. (Of course, (a) also follows
from Theorem \ref{xxsec6.7} and Lemma \ref{xxsec3.3}.)

(b) By Theorem \ref{xxsec6.7}, $kG$ is AS-regular. Since
the square of the antipode of $kG$ is the identity,
Proposition \ref{xxsec4.5}(a) shows that the Nakayama
automorphism of $kG$ is the winding automorphism induced by
the right module structure of the left integral. The
calculations in \cite[\S 3.2]{Bi1} show that this right
structure is given by $\texttt{T}_k$.

(c) This follows from Corollary \ref{xxsec5.1}.

(d) This is clear from part (b).
\end{proof}

\section*{Acknowledgments}
Both authors would like to thank Ken Goodearl, Tom Lenagan, 
Thierry Levasseur, Catharina Stroppel, Michel Van den Bergh 
and Amnon Yekutieli for useful conversations and valuable 
comments. In particular, Catharina Stroppel made many 
suggestions which have improved the exposition. The first 
author thanks Tom Hadfield for telling him about the 
``dimension drop'' phenomenon for $SL_2$. Both authors are 
supported in part by Leverhulme Research Interchange Grant 
F/00158/X (UK). The second author is supported by the National 
Science Foundation  of USA and the Royalty Research Fund 
of the University of Washington.

\providecommand{\bysame}{\leavevmode\hbox
to3em{\hrulefill}\thinspace}
\providecommand{\MR}{\relax\ifhmode\unskip\space\fi MR
}
\providecommand{\MRhref}[2]{%

\href{http://www.ams.org/mathscinet-getitem?mr=#1}{#2}
}
\providecommand{\href}[2]{#2}

\end{document}